\documentclass[a4paper,reqno]{amsart}

\usepackage{amsmath}
\usepackage{amsthm}
\usepackage{amsfonts}
\usepackage{amssymb}

\usepackage{hyperref}

\usepackage{pdfsync}

\hypersetup{
    colorlinks=false,
    pdfborder={0 0 0},
    pdfstartview={XYZ null null 1.00}
}

\numberwithin{equation}{section}

  \def\<{\langle}
  \def\>{\rangle}

  \def\ker{\mathrm{Ker}\,}
  
  \def\inv{\mathrm{Inv}\,}
  \def\inter{\mathrm{int}\,}

  \def\o{\overline}

  \def\d{ \, d }

  \def\R{\mathbb{R}}

\theoremstyle{plain}
  \newtheorem{theorem}{Theorem}[section]

  \newtheorem{proposition}[theorem]{Proposition}

  \newtheorem{corollary}[theorem]{Corollary}

\theoremstyle{definition}
  \newtheorem{definition}[theorem]{Definition}
  
  \newtheorem{remark}[theorem]{Remark}

\begin{document}

\title[On global dynamics of systems of nonlinear...]{On global dynamics of reaction--diffusion systems at resonance}

\author{Piotr Kokocki}
\address{\noindent Faculty of Mathematics and Computer Science \newline Nicolaus Copernicus University \newline Chopina 12/18, 87-100 Toru\'n, Poland}
\email{pkokocki@mat.umk.pl}


\begin{abstract}
In this paper we use the homotopy invariants methods to study the global dynamics of the reaction-diffusion systems that are at resonance at infinity. Considering degrees of the resonance for the nonlinear perturbation we establish Landesman-Lazer type conditions and use them to express the Rybakowski-Conley index of the invariant set consisting of all bounded solutions. Obtained results are applied to study the existence of solutions connecting stationary points for the system of nonlinear heat equations.
\end{abstract}
 
\maketitle

\setcounter{tocdepth}{2}

\section{Introduction}
We are concerned with the following system of autonomous differential equations 
\begin{equation}\label{sys-res}
\left\{\begin{aligned}
&\dot u(t) = -Au(t) + f(x,u(t), \nabla u(t)), && \quad (t,x)\in[0,\infty)\times\Omega,\\
&u(t) = 0, && \quad (t,x)\in[0,\infty)\times\partial\Omega,
\end{aligned}\right.
\end{equation}
where $\Omega\subset\R^{n}$ is an open bounded set with the smooth boundary and $f:\Omega\times\R^{m}\times\R^{nm}\to\mathbb{R}^{m}$ is a continuous function. Furthermore, we assume that $A$ is a linear operator given, for any $u=(u_{1},\ldots,u_{m})\in C^{2}(\o\Omega;\R^{m})$, by the formula
\begin{align*}
Au:=(A_{1}u_{1},\ldots,A_{m}u_{m}) + (\lambda_{1}u_{1},\dots,\lambda_{m}u_{m}),
\end{align*}
where $\lambda_{1},\lambda_{2},\ldots,\lambda_{m}$ are real numbers and\,\footnote{\,In the paper we use the Einstein summation convention in the definitions of differential operators.}
\begin{align*}
A_{k}u_{k} := -D_{i}(a_{k}^{ij} D_{j}u_{k})\quad\text{for} \ \  u_{k}\in C^{2}(\o\Omega)
\end{align*}
is a differential operator such that $a_{k}^{ij} = a_{k}^{ji}\in C^{1}(\o\Omega)$ for $1\le i,j\le n$ and the following inequality holds 
\begin{align*}
a_{k}^{ij}(x)\xi_{i}\xi_{j}\ge c|\xi|^{2},\quad \xi=(\xi_{1},\xi_{2},\ldots,\xi_{n})\in\R^{n}, \ x\in\Omega
\end{align*}
for some $c>0$. Throughout this paper we assume that the each of the linear operators $A_{k}$ is considered on the domain 
\begin{align*}
D(A_{k}):= \mathrm{cl}_{W^{2,p}(\Omega)}\{u_{k}\in C^{2}(\o\Omega) \ | \ u_{k}(x) = 0 \ \text{for} \ x\in\partial\Omega\},
\end{align*}
where $p\ge 2$ is an exponent that will be precisely chosen later. We are interested in the existence and topological properties of the set consisting of all bounded solutions of the system \eqref{sys-res} in the case of the {\em resonance at infinity}, that is,
\begin{align}\label{res-11}
\ker(\lambda_{k}I - A_{k})\neq\{0\} \ \ \text{and} \ \ f_{k} \text{ is a bounded map for $1\le k\le m$}.
\end{align}
It is known that under the assumption \eqref{res-11}, there are examples of the nonlinear perturbations $f$ such that the semiflow associated with the system \eqref{sys-res} does not admit bounded solutions (see Remark \ref{rem-res}). In particular, it can not have even stationary points. During last years many effort has been made to study the influence of the resonance phenomena on the existence of solutions for partial differential equations. In the fundamental paper \cite{MR0267269} the Landesman-Lazer conditions were introduced to establish the existence results for the following problem
\begin{equation}
\left\{\begin{aligned}\label{equa-stat}
&D_{i}(a^{ij} D_{j}u) + \lambda u + f_{0}(x,u) = 0,&& x\in\Omega,\\
&u = 0,&& x\in\partial\Omega,
\end{aligned}\right.
\end{equation}
where $f_{0}:\Omega\times \R\to\R$ is a bounded map, $D_{i}(a^{ij} D_{j})$ is a symmetric elliptic differential operator with the Dirichlet boundary conditions, which is considered on the space $H^{2}(\Omega)\cap H^{1}_{0}(\Omega)$, and $\lambda\in\R$ is its simple eigenvalue. Assuming that the limits $f^{\pm}_{0}(x) :=\lim_{s\to\pm\infty}f_{0}(x,s)$ exist for all $x\in\Omega$, we say that the Landesman-Lazer conditions are satisfied provided the following inequality
\begin{align}\label{lan-laz}
\int_{\Omega_{+}} f^{+}_{0}u\,dx + \int_{\Omega_{-}} f^{-}_{0}u\,dx  > 0, \ \bigg(\text{resp.} \int_{\Omega_{+}} f^{+}_{0}u\,dx + \int_{\Omega_{-}} f^{-}_{0}u\,dx<0\bigg)
\end{align}
holds for any $u\in H^{2}(\Omega)\cap H^{1}_{0}(\Omega)$ satisfying the equation  $\lambda u + D_{i}(a^{ij} D_{j}u)=0$, where $\Omega_{\pm} := \{x\in\Omega \ | \ \pm u(x)>0\}$. The results of \cite{MR0267269} were improved in \cite{MR0513090} and \cite{MR0956010} by dropping the assumption concerning the simplicity of the eigenvalue $\lambda$, whereas in \cite{MR0492839} and \cite{MR0487001} the effect of the conditions \eqref{lan-laz} on the existence of multiple equilibrium points of \eqref{equa-stat} were studied. The problem concerning the smoothness of the solutions of this equation, obtained under the Landesman-Lazer conditions, was considered in \cite{MR1466317}, where the results stating the regularity in the Besov and Triebel-Lizorkin spaces were derived. We also refer the reader to the papers \cite{MR1430505}, \cite{MR1756573}, \cite{MR1726752}, \cite{MR1689320}, where the resonance conditions were successfully adapted to finding solutions for the $p$-Laplace counterpart of the equation \eqref{equa-stat} using variational methods and linking type argument. Surprisingly, it appears that the Landesman-Lazer conditions can be applied to the study of the global dynamics of partial differential equations. To be more precise, in \cite{arieta} the results concerning the existence of global attractors for the heat equation with nonlinear boundary conditions were obtained, whereas in \cite{Kok5} the existence of orbits connecting stationary points for the parabolic partial differential equations was proved. Recently, in \cite{MR3883315} the resonance conditions \eqref{lan-laz} were used to study the existence of bifurcations from infinity for the solutions of the semilinear Schr\"odinger equation on $\R^{n}$. We refer the reader to \cite{MR2085538} and \cite{MR2609541} for the analogous bifurcation problem for $p$-Laplace and Hamilton-Jacobi-Bellman equations, respectively.
Clearly, the conditions \eqref{lan-laz} do not work if the equation \eqref{equa-stat} is at {\em strong resonance} at infinity, which means that $f(x,s)\to 0$ as $|s|\to\infty$. To handle with this case, topological and variational methods were applied in \cite{MR1041504}, \cite{MR1296115}, \cite{MR1219185}, \cite{MR713209}, \cite{MR0541873}, \cite{MR1055653}, \cite{MR0566069}, \cite{MR0612616}, \cite{MR597281} to prove the existence of solutions for the equation \eqref{equa-stat}, under various resonance assumptions imposed on the perturbation $f$. These studies were continued in \cite{MR1399539} for the $p$-Laplace version of the problem \eqref{equa-stat}. As a result, criteria were obtained on the existence of positive solutions for the $p$-Laplace equation in the terms of the sign of the limit $c:=\lim_{|s|\to\infty}sf(x,s)$, which is assumed to be independent on $x\in\Omega$. As far as we know, there is much less results in the literature concerning the resonance phenomena in case of the systems of partial differential equations. 
However, in the paper \cite{MR2280977}, the topological degree methods were applied to study the existence of periodic solutions for the following quasilinear system 
\begin{equation}\label{sys-res-per}
\left\{\begin{aligned}
&[\phi(w'(x))]' = g (x,w(x), w'(x)),  \ \ x\in(0,T), \\
&w(0) = w(T), \ \ w'(0) = w'(T), 
\end{aligned}\right.
\end{equation}
where the mapping $\phi:\R^{m}\to\R^{m}$ satisfies appropriate monotonicity and growth assumptions and $g:(0,T)\times\R^{m}\times\R^{m}\to\R^{m}$ are bounded continuous functions. It was shown that the system \eqref{sys-res-per} has a solution provided the nonlinear perturbation $g=(g_{1},\ldots,g_{m})$ satisfies the following Landesman-Lazer type conditions
\begin{equation*}
\int_{0}^{T} g_{k}^{+}(x)\,dx < 0 < \int_{0}^{T} g_{k}^{-}(x)\,dx,\quad 1\le k\le m.
\end{equation*}
In the inequalities, the function $g_{k}^{\pm}:(0,T)\to\R$, for $1\le k\le m$, is given by 
\begin{equation*}
g_{k}^{\pm}(x) := \lim_{s\to \pm\infty} g_{k}(x,u + se_{k},y),\quad x\in\Omega,
\end{equation*}
where $\{e_{j} \ | \ 1\le j\le m\}$ is the standard Euclidean basis and the limit is assumed to be uniform with respect to $u\in\mathrm{span}\{e_{j} \ | \ j\neq k\}$ and $y\in\R^{n}$. 
Motivated by the above results we intend to consider more general Landesman-Lazer type resonance conditions for the system of differential equations \eqref{sys-res}.
To this end, we assume that $1\le l\le m$ is an integer number, which is required to be fixed throughout the paper. Let $\sigma_{1},\ldots,\sigma_{m}\in[0,1]$ be given real numbers with the property that 
\begin{equation}\label{def-sig}
\check\sigma_{1}:=\min\{\sigma_{1},\ldots,\sigma_{l}\}<1,\quad \check\sigma_{2}:=\min\{\sigma_{l+1},\ldots,\sigma_{m}\}<1
\end{equation}
and let $f^{\pm}_{k}:\Omega\to\R$, for $1\le k\le m$, be continuous functions such that 
\begin{equation}\label{conv-11-aa-kk}
f_{k}^{\pm}(x) := \lim_{s\to \pm\infty} |s|^{\sigma_{k}} f_{k}(x,u+se_{k},y),\quad x\in\Omega, 
\end{equation}
where the limit \eqref{conv-11-aa-kk} is uniform with respect to $u\in \mathrm{span}\,\{e_{j} \ | \ j\neq k\}$ and $y\in\R^{nm}$. 
 We define $J_{1}$ (resp. $J_{2}$) to be the set of $1\le j\le l$ (resp. $l+1\le j\le m$) such that $\check\sigma_{1}=\sigma_{j}$ (resp. $\check\sigma_{2}=\sigma_{j}$). \label{def-l} Then we consider the following resonance conditions \\[-5pt]
$$\leqno{(LL1)_{\pm}} \left\{\begin{aligned}
&\!\sum_{k \in J_{1}}\left(\pm\!\int_{\{u_{k}>0\}}\! f_{k}^{+}(x) |u_{k}(x)|^{1-\sigma_{k}}dx \mp\int_{\{u_{k}<0\}} f_{k}^{-}(x) |u_{k}(x)|^{1-\sigma_{k}}dx \right)\!>\! 0 \\[3pt] & \text{for non-zero $(u_{1},\ldots,u_{l})\in\ker(\lambda_{1}I-A_{1})\times\ldots\times\ker(\lambda_{l}I-A_{l})$,}\end{aligned}\right.$$ 
and $$\leqno{(LL2)_{\pm}} \left\{\begin{aligned}&\!\!\sum_{k\in J_{2}}\left(\pm\!\int_{\{u_{k}>0\}} \!f_{k}^{+}(x) |u_{k}(x)|^{1-\sigma_{k}}dx \mp\!\int_{\{u_{k}<0\}} \!f_{k}^{-}(x) |u_{k}(x)|^{1-\sigma_{k}}dx\right)\!>\! 0 \\[3pt] & \text{for non-zero $(u_{l+1},\ldots,u_{m})\!\in\!\ker(\lambda_{l+1}I-A_{l+1})\!\times\!\ldots\!\times\!\ker(\lambda_{m}I-A_{m})$.}\end{aligned}\right.$$ 
Let us observe that in the assumption \eqref{conv-11-aa-kk} we use the parameter $\sigma_{k}$ to measure the strength of the resonance for the component $f_{k}$ of the nonlinear perturbation. In particular, the cases $\sigma_{k}=1$ and $\sigma_{k}=0$ correspond to the known situations that were studied for a single differential equation in \cite{MR1399539} and \cite{MR0267269}, respectively. 
To the best of our knowledge, the intermediate case $\sigma_{k}\in(0,1)$ was not considered in the literature so far. 
The main results of this paper are Theorems \ref{th-1} and \ref{th-2} that express the homotopy index of the set consisting of all bounded solutions of the system \eqref{sys-res} in the terms of the resonance conditions $(LL1)_{\pm}$ and $(LL2)_{\pm}$. 
The theorems corresponds to the earlier result of \cite{MR798176} and \cite{MR1992823}, where the index formula were obtained for the parabolic differential equation defined on a bounded domain and the whole Euclidean space $\R^{n}$, respectively, with the assumption that the resonance at infinity does not occur. 
The homotopy invariant that we use in our studies was developed in \cite{MR637695} and \cite{MR910097}, as an infinite dimensional generalization of the classical Conley index for the flows defined on finite dimensional spaces (see \cite{MR511133}, \cite{MR797044} and \cite{MR688146} for more details). The main advantage coming from this application of the homotopy invariant is that we do not require the system \eqref{sys-res} to have a gradient form, which in turn is a crucial assumption in the variational approach. Furthermore, we believe  that the conditions $(LL1)_{\pm}$ and $(LL2)_{\pm}$ are suitable to study the resonance phenomena for another systems of partial differential equations.

In the proof of Theorems \ref{th-1} and \ref{th-2} we exploit the spectral theorem for the operator $A$ to obtain a direct sum decomposition of the space $L^{p}(\Omega;\R^{m})$ into three components, among which we have the kernel of the operator $A$ and two other spaces corresponding to the positive and negative part of the spectrum of $A$. Using the homotopy invariance of the Rybakowski-Conley index we can deform the semiflow associated with the system \eqref{sys-res} to a product of semiflows defined on the spaces coming from the spectral decomposition. One of them is the $C_{0}$ semigroup generated by a restriction of the operator $-A$, while the other one is the semiflow associated with the vector field obtained by the projection of the nonlinear perturbation $f$ onto the kernel $\ker A$. Then the crucial point of our argument is to determine the contribution to the homotopy index coming from the latter semiflow, which appears to be dependent on the conditions $(LL1)_{\pm}$ and $(LL2)_{\pm}$. \\[5pt]
\noindent{\bf Outline.} The paper is organized as follows. In Section 2 we set the abstract framework to define the semiflow $\Phi$ associated with the reaction-diffusion system \eqref{sys-res} and furthermore, we recall the definition and properties of the Rybakowski-Conley homotopy index. Section 3 is devoted to the spectral decomposition of the operator $A$ on the space $L^{p}(\Omega;\R^{m})$. In Section 4 we state the main results of the paper and construction of the family of semiflows $\{\Psi^{s}\}_{s\in[0,1]}$ that will be used as the homotopy deformation of $\Phi$. In Section 5 we apply the resonance conditions $(LL1)_{\pm}$ and $(LL2)_{\pm}$ to obtain the guiding function estimates for the nonlinear perturbation $f$, whereas in Section 6 we establish a priori estimates for the bounded full solutions of the family $\{\Psi^{s}\}_{s\in[0,1]}$. Then, in Section 7, we provide the proof of Theorems \ref{th-1} and \ref{th-2} and finally, in Section 8 we provide applications of the obtained results to study the existence of solutions connecting stationary points for the system of nonlinear heat equations. \\[5pt]
\noindent {\bf Acknowledgements.} The author would like to thank the referees for their helpful suggestions and comments. 

\section{Abstract framework and homotopy index}
Let us write $X:=[L^{p}(\Omega)]^{m}$ for the real vector space equipped with the norm
\begin{align*}
\|u\|^{p}:= \sum_{k=1}^{m} \int_{\Omega} |u_{k}(x)|^{p}\,dx,\quad u=(u_{1},\ldots,u_{m})\in X
\end{align*}
and assume that $A:=(A_{1}-\lambda_{1}I)\times\ldots\times (A_{m}-\lambda_{m}I)$ is the product operator defined on the space $X$. It is known (see e.g. \cite{MR1778284}, \cite{MR610244}, \cite{Pazy}, \cite{MR500580}) that the operator $A$ is sectorial, that is, there are $\gamma\in(0,\pi/2)$, $C_{1}\ge 1$ and $a\in\mathbb{R}$, such that the sector $$\Sigma_{a,\gamma}:=\{\lambda\in\mathbb{C} \ | \ \gamma \le |\mathrm{arg} \, (\lambda - a)| \le \pi, \ \lambda\neq a\}$$ is contained in the resolvent set $\rho(A)$ of the operator $A$ and the inequality holds $$\|(\lambda I - A)^{-1}\| \le C_{1}/|\lambda - a|, \quad \lambda\in \Sigma_{a,\gamma}.$$ Furthermore $-A$ is a generator of a compact analytic semigroup $\{S_{A}(t)\}_{t\ge 0}$ of bounded linear operators on $X$. Let us observe that from \cite[Theorem 16.7.2]{H-F} we have the following useful kernel relation
\begin{equation}\label{eq-ker}
\ker A = \ker(I - S_{A}(t)),\quad t>0.
\end{equation}
If we write $\delta:=1+\max\{\lambda_{k} \ | \ 1\le k\le m\}$, then the operator $A_{\delta}:=A+\delta I$ is positively defined, that is, $\mathrm{Re}\,\mu > 0$ for any element $\mu$ from the spectrum $\sigma(A+\delta I)$. Hence, given $\alpha \ge 0$, we can define {\em the fractional space} $X^\alpha$ as the domain of the fractional power $(\delta I + A)^{\alpha}$ (see \cite[Section 2.6]{Pazy}), endowed with the graph norm 
\begin{equation*}
\|u\|_\alpha := \|(A+\delta I)^{\alpha} u\|,\quad u\in X^{\alpha}.
\end{equation*}
It is known that $X^{\alpha}$ is a Banach space, continuously embedded in $X$, that is, there is a constant $C_{2}>0$ such that the following inequality holds
\begin{align}\label{cont-emb}
\|u\|\le C_{2}\|u\|_{\alpha},\quad u\in X^{\alpha}.
\end{align}
From now on we assume additionally that 
\begin{equation}\label{alpha-p}
\alpha\in(3/4,1), \quad p\ge 2n
\end{equation}
and furthermore we require that, for any $1\le k\le m$, the continuous map $f_{k}:\Omega\times\R^{m}\times\R^{mn}\to \R$ satisfies the following conditions: \\[3pt] 
\noindent\makebox[22pt][l]{$(F1)$} \parbox[t][][t]{118mm}{given $R>0$, there exists a constant $L_{R} > 0$ such that
\begin{align*}
|f_{k}(x,s_1,y_1) - f_{k}(x,s_2,y_2)| & \le L_{R}\,(|s_1 - s_2| + |y_1 - y_2|),
\end{align*}
for $x\in\Omega$, $s_1,s_2\in\R^{m}$ and $y_1,y_2\in\mathbb{R}^{nm}$ with $|s_{1}|,|s_{2}|,|y_{1}|,|y_{2}|\le R$;}\\[5pt]
\noindent\makebox[22pt][l]{$(F2)$} \parbox[t][][t]{118mm}{there exists a constant $C_{3}>0$ such that, for any $1\le k\le m$, we have $$|f_{k}(x,s,y)| \le C_{3},  \quad x\in\Omega, \ s\in\R^{m}, \ y\in\mathbb{R}^{mn}.$$}\\[2pt]
In view of the assumption \eqref{alpha-p}, we can check that $2\alpha - n/p>1$, which by \cite[Theorem 1.6.1]{MR610244}, gives the inclusion $X^{\alpha}\subset C^{1}(\o\Omega)$ together with the inequality
\begin{equation*}
\|u\|_{C^{1}(\o\Omega)} \le C_{4}\|u\|_{\alpha},\quad u\in X^{\alpha},
\end{equation*}
where $C_{4}>0$ is a constant. Therefore the map $F:X^{\alpha}\to X$, given by the formula 
\begin{align*}
F(u):=(f_{1}(x, u,\nabla u), \ldots, f_{m}(x, u,\nabla u)),\quad u\in X^{\alpha}
\end{align*}
is well-defined and straightforward calculations show that it is bounded and satisfies the Lipschitz condition on the bounded subsets of $X^{\alpha}$. Consequently the system \eqref{sys-res} can be written in the following abstract form 
\begin{equation}\label{row-a-f}
\dot u(t) = - A u(t) + F (u(t)), \quad  t > 0.
\end{equation}
\begin{definition}(see \cite[page 57]{MR1778284}, \cite[page 53]{MR610244})
Given the interval $I\subset\R$, we say that the function $u:I\to X^{\alpha}$ is a \emph{mild solution} of the equation \eqref{row-a-f}, provided
\begin{equation*}
u(t) = S_A(t - t')u(t') + \int_{t'}^t S_A(t - \tau)F(u(\tau)) \,d \tau \ \ \text{for} \ \ t,t'\in I \ \ \text{with} \ \ t'<t.
\end{equation*}
\end{definition}
From \cite[Theorem 3.3.3]{MR610244}, \cite[Corollary 3.3.5]{MR610244} and Remark \ref{rem-lipsch-bound} it follows that, for any $u_{0}\in X^{\alpha}$, equation \eqref{row-a-f} admits a unique global mild solution $u(\,\cdot\,;u_{0}):[0,+\infty)\to X^\alpha$ starting at $u_{0}$. Hence we can define a semiflow $\Phi:[0,+\infty)\times X^\alpha \to X^\alpha$ associated with the equation \eqref{row-a-f} by the following formula $$\Phi(t,u_{0}):= u(t;u_{0}),\quad t\ge 0, \  u_{0} \in X^\alpha.$$ From \cite[Proposition 2.3.2]{MR1778284}, we infer that the semiflow is continuous, that is, for any sequence $( u_{n})$ in $X^{\alpha}$ such that $ u_{n}\to u_{0}$ as $n\to\infty$, we have
\begin{align*}
\Phi(t; u_{n})\to \Phi(t;  u_{0})\quad \text{for} \ \ t\ge 0, \ \ \text{as} \ \ n\to\infty
\end{align*}
and the convergence is uniform for the time $t$ from bounded subsets of the interval $[0,+\infty)$. Furthermore it is known that the operator $A$ has compact resolvents and therefore, we can use \cite[Theorem 3.2.1]{MR1778284} to deduce that any bounded set $M\subset X^\alpha$ is {\em admissible} with respect to $\Phi$, which means that, for every sequences $( u_n)$ in $X^\alpha$ and $(t_n)$ in $[0,+\infty)$, if  $t_n \to +\infty$ as $n\to \infty$ and $$\Phi([0,t_n]\times \{ u_n\})\subset M,\quad  n\ge 1,$$ then the set $\{\Phi(t_n, u_n) \ | \ n\ge 1\}$ is relatively compact in the space $X^\alpha$. 
\begin{definition}(see e.g. \cite[page 3]{MR910097})
Assume that $u:[-\delta_1, \delta_2) \to X^\alpha$, where $\delta_1 \ge 0$ and $\delta_2 > 0$, is a continuous map. We say that the map $u$ is a \emph{solution} of the semiflow $\Phi$, provided $$\Phi(t,u(s)) = u(t+s) \ \ \text{for} \ \ t\ge 0 \ \ \text{and} \ \ s\in [-\delta_1,\delta_2) \ \ \text{such that} \ t + s\in[-\delta_1,\delta_2).$$ In particular, if the map $u$ is defined on the whole real line, then $u$ is called a {\em full solution} of the semiflow $\Phi$. \hfill $\square$
\end{definition}
We recall that the set $K\subset X^{\alpha}$ \emph{invariant} with respect to $\Phi$ provided, for every $u_{0}\in K$ there is a full solution $u$ of the semiflow $\Phi$ such that $u(0) =  u_{0}$ and $u(\R) \subset K$. Therefore, given $M\subset X^\alpha$, we define its maximal invariant subset $\inv M = \inv(M,\Phi)$ as the set of points $u_{0}\in M$ with the property that there is a full solution $u$ of the semiflow $\Phi$ such that $u(0) =  u_{0}$ and $u(\R)\subset M$. In particular, we call $K$ an \emph{isolated invariant set}, if there is a closed set $M\subset X^\alpha$ such that 
\begin{equation*}
K=\inv M\subset\inter M.
\end{equation*}
Then we say that $M$ is an \emph{isolating neighborhood} for $K$. 
\begin{definition}(see \cite[page 9]{MR910097})
Assume that $B\subset X^\alpha$ is a closed set and let $u_{0}\in\partial B$. We say that $u_{0}$ is \emph{a strict egress point} (resp. \emph{strict ingress point}, resp. \emph{bounce off point}), if for any solution $u:[-\delta_1, \delta_2) \to X^\alpha$, where $\delta_1 \ge 0$ and $\delta_2 > 0$, of the semiflow $\Phi$ such that $u(0) =  u_{0}$ the following conditions are satisfied:
\begin{enumerate}
\item[(a)] there is $\varepsilon_2\in(0,\delta_2]$ such that $u(t)\notin B$ (resp. $u(t)\in \inter B$, resp. $u(t)\notin B$) for $t\in(0,\varepsilon_2]$;
\item[(b)] if $\delta_1 > 0$ then there is $\varepsilon_1\in(0,\delta_1)$ such that $u(t)\in \inter B$ (resp. $u(t)\notin B$, resp. $u(t)\notin B$) for $t\in[-\varepsilon_1, 0)$.
\end{enumerate}
Then we write $B^e$ (resp. $B^i$, resp. $B^b$) for the set of strict egress points (resp. strict ingress points, resp. strict bounce off points) and furthermore, we set $B^-:=B^e \cup B^b$. We say that the closed set $B\subset X^\alpha$ is an \emph{isolating block}, provided $\partial B = B^e\cup B^i\cup B^e$ and $B^-$ is a closed set in $X^{\alpha}$. \hfill $\square$
\end{definition}
Throughout the paper, we write $[Y_{0},y_{0}]$ for the homotopy type of the pointed topological space $(Y_{0},y_{0})$. In particular, if $Y_{0}=\{y_{0}\}$ then we say that the homotopy type $[\{y_{0}\},y_{0}]$ is trivial and we denote it by $\o 0$. Furthermore, given $k\ge 0$, we set $\Sigma^{k}:=[S^{k},s_{0}]$, where $S^{k}$ is $k$-dimensional unit sphere and $s_{0}\in S^{k}$ is an arbitrary point. Given two pointed spaces 
$(Y_{0},y_{0})$ and $(Y_{1},y_{1})$ we define their {\em smash product} $(Y_{0},y_{0})\wedge(Y_{1},y_{1})$ as a pair $(Y_{2},y_{2})$, where $Y_{2}:=(Y_{0}\times Y_{1})/(Y_{0}\times\{y_{1}\}\cup \{y_{0}\}\times Y_{1})$ is the quotient topological space and $y_{2}:=[Y_{0}\times\{y_{1}\}\cup \{y_{0}\times Y_{1}\}]$ is the distinguished point. It is well-known that the homotopy type of the smash product of two pointed spaces actually depends on the homotopy type of each of them. Therefore, we have the well-defined operation $[Y_{0},y_{0}]\wedge[Y_{1},y_{1}]:=[(Y_{0},y_{0})\wedge(Y_{1},y_{1})]$ (see e.g. \cite{bredon}, \cite{hat}, \cite{MR910097} and \cite{spanier} for more details).\\[5pt]
\indent From \cite[Theorem I.5.1]{MR910097} we know that, for any isolated invariant set $K$, which admits an admissible isolating neighborhood, we can construct an isolating block $B$ such that $K=\inv B$. Then we define the {\em homotopy index} of $K$ as 
\begin{equation*}
h(\Phi, K):=\left\{\begin{aligned}&[B/B^-, [B^-]], && \text{ if } B^-\neq\emptyset, \\ & [B\,\dot\cup\,\{c\}, c], && \text{ if } B^- = \emptyset,\end{aligned}\right.
\end{equation*}
where $B/B^-$ is the quotient topological space and $B\,\dot\cup\,\{c\}$ is a disjoint sum of $B$ and the one point space $\{c\}$. 
It is known that the homotopy index is independent on the choice of isolating block $B$ for the set $K$ and has the following properties.\\[5pt]
\makebox[8mm][l]{(H1)}\parbox[t]{118mm}{If $M\subset X^{\alpha}$ is an admissible isolating neighborhood and the homotopy index of $K:=\inv M$ is nontrivial, then the set $K$ is non-empty.}\\[5pt]
\makebox[8mm][l]{(H2)}\parbox[t]{118mm}{Let $\varphi_{j}:[0,+\infty)\times Z_{j}\to Z_{j}$, for $j=1,2$, be semiflows defined on the closed components of the direct sum decomposition $X^{\alpha} = Z_{1}\oplus Z_{2}$. Assume that, for any $j=1,2$, the set $M_{j}\subset Z_{j}$ is an admissible isolating neighborhood for $K_{j}:=\inv(M_{j},\varphi_{j})$. Then the set $M_{1}\oplus M_{2}$ is an admissible isolating neighborhood with respect to the product semiflow $\varphi_{1}\oplus\varphi_{2}$ and 
$$h(\varphi_{1}\oplus\varphi_{2},K)=h(\varphi_{1},K_1)\wedge h(\varphi_{2},K_2),$$ where 
$K:=\inv(M_{1}\oplus M_{2},\varphi_{1}\oplus\varphi_{2})$.}\\[5pt]
\makebox[8mm][l]{(H3)}\parbox[t]{118mm}{Assume that the closed set $M\subset X^{\alpha}$ is admissible with respect to the family of semiflows 
$\{\Psi^{s}\}_{s\in[0,1]}$, that is, for every sequences $(s_n)$ in $[0,1]$, $(u_n)$ in $X^\alpha$ and $(t_n)$ in $[0,+\infty)$, if $t_n \to +\infty$ as $n\to \infty$ and $$\Psi^{s_n}([0,t_n]\times\{u_n\})\subset M,\quad  n\ge 1,$$ then the set $\{\Psi^{s_n}(t_n, u_n) \ | \ n\ge 1\}$ is relatively compact in $X^\alpha$. If the set $M$ is an isolating neighborhood of $K_s:=\inv(\Psi^{s}, M)$ for all $s\in[0,1]$, then 
$$h(\Psi^{0}, K_{0}) = h(\Psi^{1}, K_{1}).$$}\\[5pt]
Let us assume that $u$ is a full solution of the semiflow $\Phi$. We define the limit set $\alpha(u)$ (resp. $\omega(u)$) to be the collection of points $u'\in X^{\alpha}$ such that $u(t_{n})\to  u'$ for some $t_{n}\to-\infty$ (resp. $t_{n}\to+\infty$). The following proposition (see \cite[Theorem 11.5]{MR910097}) provides a tool to study the existence of solutions connecting invariant sets by the use of the homotopy index.
\begin{proposition}\label{prop-irred}
Let us assume that $K$ and $K_{0}$ are isolated invariant sets such that each of them admits an admissible isolating neighborhood and furthermore
\begin{align*}
h(\varphi, K_{0})\neq \o 0 \ \ \text{and} \ \ h(\varphi, K) = \Sigma^{k} \ \ \text{for some} \ \ k\ge 0.
\end{align*}
If $h(\varphi, K_{0})\neq h(\varphi, K)$ then there is a full solution $u$ of the semiflow $\Phi$ such that either $\alpha(u)\subset K_{0}$ or $\omega(u)\subset K_{0}$.
\end{proposition}

\section{Spectral decomposition of the product operator}
In this section we assume that $\lambda:=(\lambda_{1},\ldots,\lambda_{m})$ is a vector consisting of real numbers. Let us observe that, for any $1\le k\le m$, the spectrum $\sigma(A_{k})$ consists of a bounded below sequence of real eigenvalues, which is either finite or divergent to infinity. Since the respective eigenspaces of the operator $A_{k}$ are finite dimensional, we can define the number $d_{k}(\lambda):=0$ if $\lambda_{k}$ is the first eigenvalue of $A_{k}$ and otherwise
\begin{gather*}
d_{k}(\lambda):=\sum_{\substack{\nu<\lambda_{k}}} \dim\ker(\nu I - A_{k}),
\end{gather*}
where in the above summation the parameter $\lambda$ ranges over the set $\sigma_{p}(A_{k})$ of eigenvalues of the operator $A_{k}$. Then we put 
\begin{gather}
d_{\infty}(\lambda):=d_{1}(\lambda)+\ldots+d_{m}(\lambda).
\end{gather}
Let us consider the auxiliary product operator $\widehat A:=(\widehat A_{1}-\lambda_{1}I)\times\ldots\times (\widehat A_{m}-\lambda_{m}I)$, defined on the space $\widehat X:=L^{2}(\Omega;\R^{m})$, where, for any $1\le k\le m$, we assume that $\widehat A_{k}$ is a linear operator on $L^{2}(\Omega)$ given by
\begin{equation}\label{sect-op-2}
\left\{\begin{aligned}
D(\widehat A_{k})&:= \mathrm{cl}_{W^{2,2}(\Omega)}\{u_{k}\in C^{2}(\o\Omega) \ | \  u_{k}(x) = 0\text{ for }x\in\partial\Omega\},\\
\widehat  A_{k}u_{k}&:= D_{i}(a^{ij}_{k} D_{j}u_{k}), \quad u_{k}\in D(\widehat  A_{k}).
\end{aligned}\right.
\end{equation}
\begin{remark}\label{sect-op-2-a}
We claim that $\sigma_{p}(A)=\sigma_{p}(\widehat{A})\subset\R$. Indeed, since $\widehat A$ is a symmetric operator and $A\subset \widehat A$, its spectrum consists of real eigenvalues. Furthermore, the fact that $A\subset \widehat A$ implies $\sigma_{p}(A)\subset \sigma_{p}(\widehat A)$. To check the opposite inclusion let us observe that, by the regularity properties of the elliptic operators (see e.g. \cite{MR500580}), the eigenvalues of the operator $\widehat A$ can be considered as smooth functions on the set $\o\Omega$. Therefore, if $\mu \in\sigma_{p}(\widehat A)$ and $u\in \mathrm{Ker}\,(\mu I - \widehat A)$, then $u\in D(A)$ and $(\mu I - A)u = (\mu I - \widehat A)u=0$. Hence $\mu\in\sigma_{p}(A)$ and $\sigma_{p}(\widehat A)\subset \sigma_{p}(A)$ as desired. \hfill $\square$
\end{remark}
In view of Remark \ref{sect-op-2-a} and the fact that the operator $A$ has compact resolvents, the spectrum $\sigma(A)$ can be represented as the sequence of real isolated eigenvalues $(\mu_{k})_{k\ge 1}$, which is finite or $\mu_{k}\to\infty$ as $k\to\infty$ (see e.g. \cite{MR1778284}, \cite{MR500580}). On the other hand, the resonance assumption \eqref{z-zero}, gives $\mathrm{Ker}\,A\neq \{0\}$, and hence we can choose $r\ge 1$ such that 
\begin{align*}
\mu_{1}<\ldots<\mu_{r-1}<\mu_{r} = 0\quad\text{and}\quad 0<\mu_{k} < \mu_{k+1},\quad k\ge r+1.
\end{align*}
By the spectral theorem for the symmetric operators with compact resolvents (see \cite[Theorem 1.5.2]{MR610244}) we obtain the direct sum decomposition $\widehat X = \widehat X_{-}\oplus\widehat X_{0}\oplus \widehat X_{+}$ on the closed and mutually orthogonal subspaces of $\widehat X$, where
\begin{align}\label{ww11}
\widehat X_{0} = \mathrm{Ker}\,(\mu_{r}I - \widehat A)\quad\text{and}\quad \widehat X_{-} = \mathrm{Ker}\,(\mu_{1}I - \widehat A)\oplus\ldots\oplus \mathrm{Ker}\,(\mu_{r-1}I - \widehat A).
\end{align}
Furthermore we have the following inclusions
\begin{equation}\label{incl-2233}
\widehat A(\widehat X_{0})\subset\widehat X_{0},\quad \widehat A(\widehat X_{-})\subset\widehat X_{-},\quad \widehat A(D(A)\cap \widehat X_{+})\subset\widehat X_{+}
\end{equation}
and, if $\widehat A_{k}$ is a part of the operator $\widehat A$ in the space $\widehat X_{k}$, then
\begin{align}\label{spec-eq}
\sigma(\widehat A_{-}) = \{\mu_{1},\ldots,\mu_{r-1}\}\quad\text{and}\quad \sigma(\widehat A_{+}) = \{\mu_{k} \ | \ k\ge r+1\}.
\end{align}
Let us define $X_{k}:=X\cap \widehat X_{k}$ for $k\in\{0,-,+\}$ and write
\begin{equation}\label{def-nn}
\begin{gathered}
N_{1}:= \mathrm{Ker}\,(A_{1}-\lambda_{1}I)\times\ldots\times\mathrm{Ker}\,(A_{l}-\lambda_{l}I),\\
N_{2}:= \mathrm{Ker}\,(A_{l+1}-\lambda_{l+1}I)\times\ldots\times\mathrm{Ker}\,(A_{m}-\lambda_{m}I),
\end{gathered}
\end{equation}
where we recall that $1\le l\le m$ is the fixed integer number used in the definition of the resonance conditions $(LL1)_{\pm}$ and $(LL2)_{\pm}$ (see page \pageref{def-l}).
\begin{remark}\label{rem-dim-eq}
We claim that
\begin{equation*}
\mathrm{dim}\, X_{-} = d_{\infty}(\lambda).
\end{equation*}
Indeed, as an immediate consequence of \eqref{ww11}, we obtain 
\begin{align}\label{eq-hhpp}
X_{-} = \widehat X_{-} = \mathrm{Ker}\,(\mu_{1}I - \widehat A)\oplus\ldots\oplus \mathrm{Ker}\,(\mu_{r-1}I - \widehat A).
\end{align}
Furthermore, for any $1\le k\le r-1$, we have
\begin{align}\label{eq-hhpp-22}
\mathrm{Ker}\,(\mu_{k}I - \widehat A) = \mathrm{Ker}\,((\mu_{k}+\lambda_{1})I - \widehat A_{-})\times\ldots\times\mathrm{Ker}\,((\mu_{k}+\lambda_{m})I - \widehat A_{m}).
\end{align}
Let us observe that, given $1\le j\le m$, we have
\begin{gather*}
\{\mu_{k} + \lambda_{j} \ | \ 1\le k\le r-1\text{ and }\mathrm{Ker}\,((\mu_{k} + \lambda_{j})I - \widehat A_{j}) \neq \{0\}\} \\
 = \{\lambda<\lambda_{j} \ | \ \mathrm{Ker}\,(\lambda I - \widehat A_{j})\neq\{0\}\},
\end{gather*}
which together with \eqref{eq-hhpp} and \eqref{eq-hhpp-22} give
\begin{align*}
\mathrm{dim}\, X_{-} &= \sum_{k=1}^{r-1}\mathrm{dim}\,\mathrm{Ker}\,(\mu_{k}I - \widehat A) = \sum_{k=1}^{r-1}\sum_{j=1}^{m}\mathrm{dim}\,\mathrm{Ker}\,((\mu_{k} + \lambda_{j})I - \widehat A_{j})\\
& = \sum_{j=1}^{m}\sum_{\nu<\lambda_{j}}\mathrm{dim}\,\mathrm{Ker}\,(\nu I - \widehat A_{j}) = \sum_{j=1}^{m} d_{j}(\lambda) = d_{\infty}(\lambda),
\end{align*}
as claimed. \hfill $\square$ 
\end{remark}
Let us observe that $X_{k} = \widehat X_{k}\subset D(A)$ for $k\in\{0,-\}$ and the space $X$ can be represented as the direct sum $X=N_{1}\oplus N_{2}\oplus X_{-}\oplus X_{+}$, where the component spaces are closed in $X$ and mutually orthogonal in $\widehat X = L^{2}(\Omega;\R^{m})$. For any $k=1,2$ and $l\in\{-,+\}$, we denote by $P_{k}$ and $Q_{l}$ the projections on the spaces $N_{k}$ and $X_{l}$, respectively, that are determined by this decomposition. We also write $Q_{0} := P_{1}+P_{2}$. Then, by the continuity of the inclusion $X^\alpha\subset X$, we obtain $$X^\alpha= N_{1}\oplus N_{2}\oplus X^\alpha_{-}\oplus X^\alpha_{+},$$ where $X^\alpha_{k}:=X^\alpha\cap X_{k}$ for $k\in\{-,+\}$, are closed subspaces of $X^{\alpha}$. Hence, the linear operator $Q_{k}$ can be restricted to the bounded map $Q_{k}:X^{\alpha}\to X^{\alpha}$ for $k\in\{-,+\}$.
\begin{remark}\label{rem-aa-1}
Given $k\in\{-,+\}$, we denote by $A_{k}$ the part of the operator $A$ in the space $X_{k}$. We claim that $\rho(A)\subset\rho(A_{k})$ and 
\begin{equation}\label{eq-4455}
(\mu I-A)^{-1}v_{k} = (\mu I-A_{k})^{-1}v_{k},\quad v_{k}\in X_{k}.
\end{equation}
Indeed, let us assume that $k\in\{-,+\}$ is fixed and take $\mu\in\rho(A)$. Then we have $\mathrm{Ker}\,(\mu I - A_{k})=\{0\}$ because $A_{k}\subset A$. On the other hand, if we take arbitrary $v\in X_{k}$ then $(\mu I-A)u = v$ for some $u\in D(A)$. Let us write $u=u_{0}+u_{1}+u_{2}$, where $u_{i}\in X_{i}$ for $i\in\{0,1,2\}$. Since $u,u_{0},u_{1}\in D(A)$, it follows also that $u_{2}\in D(A)$ and
\begin{align*}
v = (\mu I- A)u = (\mu I-\widehat A)u = (\mu I-\widehat A)u_{0} + (\mu I-\widehat A)u_{1}+(\mu I-\widehat A)u_{2}.
\end{align*}
Combining this with \eqref{incl-2233} and the fact that $v\in X_{k}$, we have $v = (\mu I-\widehat A)u_{k} = (\mu I-A)u_{k}$, which gives $u_{k}\in D(A_{k})$ and $(\mu I-A_{k})u_{k} = v$. Hence the operator $\mu I-A_{k}$ is invertible on $X_{k}$ and \eqref{eq-4455} holds. Observe that the operator $A_{k}$ is closed as a part of the closed operator $A$ in the closed subspace $X_{k}\subset X$. Consequently the inverse operator $(\mu I-A_{k})^{-1}$ is bounded on $X_{k}$ and $\mu\in\rho(A_{k})$ as claimed. \hfill $\square$
\end{remark}
Observe that an argument similar to that in Remark \ref{sect-op-2-a} shows that $\sigma_{p}(A_{+})=\sigma_{p}(\widehat A_{+})$. Since the operator $A$ has compact resolvents, by Remark \ref{rem-aa-1} the operator $A_{+}$ also has the property. Consequently, we have $\sigma(A_{+}) =\sigma_{p}(\widehat A_{+})=\sigma_{p}(\widehat A_{+})$, which together with \eqref{spec-eq} and the fact that $A_{-}=\widehat A_{-}$ provide
\begin{align}\label{spect-11}
\sigma(A_{-}) = \{\mu_{1},\ldots,\mu_{r-1}\}\quad\text{and}\quad \sigma(A_{+}) = \{\mu_{j} \ | \ j\ge r+1\}.
\end{align}
Let us take arbitrary $\mu\in\rho(A)$ and observe that from the definition of the operator $A$ and equality \eqref{eq-4455}, we have 
\begin{align*}
&Q_{k}(\mu I - A)^{-1}u = (\mu I - A)^{-1}Q_{k}u = (\mu I - A_{k})^{-1}Q_{k}u, && \hspace{-15pt}u\in X,\\
&P_{k}(\mu I - A)^{-1}u = (\mu I - A)^{-1}P_{k}u = (\mu I - A_{k})^{-1}P_{k}u, && \hspace{-15pt}u\in X,
\end{align*}
which by the Euler formula for the $C_{0}$ semigroups (see \cite[Theorem 8.3]{Pazy}) yields
\begin{align}\label{semig-2}
&S_{A_{k}}(t)Q_{k}u = S_{A}(t)Q_{k}u = Q_{k}S_{A}(t)u, && \hspace{-40pt} t\ge 0, \ u\in X, \\ \label{semig-2bb}
&S_{A}(t)P_{k}u = P_{k}S_{A}(t)u, && \hspace{-40pt} t\ge 0, \ u\in X.
\end{align}
\begin{remark}
The semigroup $\{S_{A}(t)\}_{t\ge 0}$ extends on the space $X_{-}$ to a $C_{0}$ group of bounded linear operators and there are constants $c,C_{5}>0$ such that 
\begin{align}\label{bb-nn-2}
& \|S_A(t)u\| \le C_{5} e^{c t}\|u\|, && \hspace{-50pt} t\le 0, \ u\in X_{-}, \\ \label{bb-nn-3}
& \|S_A(t)u\| \le C_{5} e^{-ct} \|u\|, && \hspace{-50pt}  t \ge 0, \ u\in X_{+}, \\ \label{bb-nn-1}
& \|S_A(t)u\|_{\alpha} \le C_{5} e^{-ct} t^{-\alpha} \|u\|, && \hspace{-50pt}  t > 0, \ u\in X_{+}.
\end{align}
Indeed, observe that the equality \eqref{semig-2} implies that 
\begin{equation}\label{eq-11-pp}
S_{A}(t)u = S_{A_{-}}(t)u,\quad t\ge 0, \quad u\in X_{-}.
\end{equation}
Since $X_{-}$ is a finite dimensional space, the operator $A_{-}$ generates the $C_{0}$ group of bounded linear operators on $X_{-}$, which by the equality \eqref{eq-11-pp} is the desired extension of $\{S_{A}(t)\}_{t\ge 0}$ on $X_{-}$. Furthermore, by \eqref{spect-11}, the operator $-A_{-}$ is positively definite which implies that
\begin{align*}
\|S_A(t)u\| = \|S_{A_{-}}(t)u\| \le C e^{c t}\|u\|, \quad t\le 0, \ u\in X_{-},
\end{align*}
where $c,C>0$ are constants and consequently the inequality \eqref{bb-nn-2} follows. As a direct consequence of Remark \ref{rem-aa-1} and \eqref{spect-11} we find that $A_{+}$ is positively defined sectorial operator on $X_{+}$, which by \cite[Theorem 6.13]{Pazy} allows us to modify the constants $c,C>0$ if necessary, to obtain the following estimates
\begin{align}\label{bb-nn-3a}
& \|S_{A_{+}}(t)u\| \le C e^{-ct} \|u\|, && \hspace{-40pt}  t \ge 0, \ u\in X_{+}, \\ \label{nn-ee11a}
& \|A_{+}^{\alpha}S_{A_{+}}(t)u\| \le C e^{-ct} t^{-\alpha} \|u\|, && \hspace{-40pt}  t > 0, \ u\in X_{+}.
\end{align}
On the other hand, \eqref{semig-2} implies that $S_{A}(t)u = S_{A_{+}}(t)u$ for $t\ge 0$ and $u\in X_{+}$, which together with the inequality \eqref{bb-nn-3a} gives
\begin{align*}
\|S_A(t)u\| = \|S_{A_{+}}(t)u\| \le C e^{-ct}\|u\|, \quad t\ge 0, \ u\in X_{+},
\end{align*}
and \eqref{bb-nn-3} follows. Observe that, by \cite[Theorem 1.4.6]{MR610244} we have the equality of domains $D(A_{+}^{\alpha}) = D((\delta I+A_{+})^{\alpha})$ and the equivalence of the corresponding norms 
\begin{align}\label{eq-norm-1}
C''\|A_{+}^{\alpha}u\| \le \|(\delta I+A_{+})^{\alpha}u\| \le C'\|A_{+}^{\alpha}u\|,\quad u\in D(A_{+}^{\alpha}),
\end{align}
where $C',C''>0$. On the other hand, by the definition of the fractional power of the positive sectorial operator $ (\delta I+A_{+})^{\alpha}\subset  (\delta I+A)^{\alpha}$, which together with \eqref{nn-ee11a} and \eqref{eq-norm-1}, for any $t>0$ and $u\in X_{+}$ yields
\begin{align*}
&\|S_{A}(t)u\|_{\alpha} = \|(\delta I+A)^{\alpha} S_{A}(t)u\| = \|(\delta I+A)^{\alpha} S_{A_{+}}(t)u\| \\
&\qquad = \|(\delta I+A_{+})^{\alpha} S_{A_{+}}(t)u\| \le C'\|A_{+}^{\alpha} S_{A_{+}}(t)u\| \le CC' e^{-ct} t^{-\alpha} \|u\|,
\end{align*}
and the estimate \eqref{bb-nn-1} follows. \hfill $\square$
\end{remark}

\section{Statement of the main results}
Given $\lambda=(\lambda_{1},\ldots,\lambda_{m})$ we impose the following standing resonance assumption
\begin{align}\label{z-zero}
\mathrm{Ker}\,(\lambda_{k} I - A_{k})\neq \{0\},\quad 1\le k\le m
\end{align}
and define the following numbers
\begin{gather}
n_{1}(\lambda):= \sum_{i=1}^{l} \dim\ker(\lambda_{i} I - A_{i}),\quad n_{2}(\lambda):= \sum_{i=l+1}^{m} \dim\ker(\lambda_{i} I - A_{i}).
\end{gather}
The main results of this paper are the following theorems. 
\begin{theorem}\label{th-1}
Suppose that $\{h_{k}\}_{k=1}^{m}$ is a family of $L^{2}(\Omega)$ functions such that 
$$\leqno{(C1)_{\pm}} \quad 
\left\{\begin{gathered}
\text{for any $1\le k\le l$ the following inequality holds}\\ 
\pm f_{k}(x,u,y) |u_{k}|^{\sigma_{k}} \mathrm{sgn}\,u_{k} \ge h_{k}(x), \ \ x\in\Omega, \ u\in\R^{m}, \ y\in\R^{nm},
\end{gathered}\right.$$
and furthermore
$$\leqno{(C2)_{\pm}}  \quad
\left\{\begin{gathered}
\text{for any $l+1\le k\le m$ the following inequality holds}\\ 
\pm f_{k}(x,u,y) |u_{k}|^{\sigma_{k}} \mathrm{sgn}\,u_{k} \ge h_{k}(x), \ \ x\in\Omega, \ u\in\R^{m}, \ y\in\R^{nm}.
\end{gathered}\right.$$
If the conditions $(LL1)_{\pm}$ and $(LL2)_{\pm}$ are satisfied, then the set $K_{\infty}$ consisting of all bounded full solutions of the semiflow $\Phi$ admits an admissible isolating neighborhood. Furthermore the homotopy index of $K_{\infty}$ is given by  
\begin{equation}\label{form-conley-1}
\hspace{5pt}h(\Phi, K_{\infty}) = \Sigma^{d_{\infty}(\lambda) + n_{1}(\lambda)+n_{2}(\lambda)},
\end{equation}
if the conditions $(C1)_{+}$, $(C2)_{+}$, $(LL1)_{+}$, $(LL2)_{+}$ hold and
\begin{equation}\label{form-conley-1b}
\hspace{15pt}h(\Phi, K_{\infty}) = \Sigma^{d_{\infty}(\lambda)}, 
\end{equation}
if the conditions $(C1)_{-}$, $(C2)_{-}$, $(LL1)_{-}$, $(LL2)_{-}$ are satisfied. 
\end{theorem}
In the subsequent result we prove analogous homotopy index formula in the case of the resonance conditions $(LL1)_{\pm}$ and $(LL2)_{\mp}$.
\begin{theorem}\label{th-2}
Suppose that $\{h_{k}\}_{k=1}^{m}$ is a family of $L^{2}(\Omega)$ functions such that the inequalities $(C1)_{\pm}$ and $(C2)_{\mp}$ hold. 
If conditions $(LL1)_{\pm}$ and $(LL2)_{\mp}$ are satisfied, then the set $K_{\infty}$ consisting of all bounded full solutions of the semiflow $\Phi$ admits an admissible isolating neighborhood and the homotopy index of $K_{\infty}$ is given by  
\begin{equation}\label{form-conley-2}
\hspace{15pt}h(\Phi, K_{\infty}) = \Sigma^{d_{\infty}(\lambda) + n_{1}(\lambda)},
\end{equation}
if the conditions $(C1)_{+}$, $(C2)_{-}$, $(LL1)_{+}$, $(LL2)_{-}$ hold and 
\begin{equation}\label{form-conley-2b}
\hspace{15pt}h(\Phi, K_{\infty}) = \Sigma^{d_{\infty}(\lambda) + n_{2}(\lambda)}, 
\end{equation}
if the conditions $(C1)_{-}$, $(C2)_{+}$, $(LL1)_{-}$, $(LL2)_{+}$ are satisfied. 
\end{theorem}
\begin{remark}
According to the formulas \eqref{def-sig}, we intend to prove Theorems \ref{th-1} and \ref{th-2} with the standing assumption that $0\le \check\sigma_{1},\check\sigma_{2}<1$. The calculation of the Rybakowski-Conley index of the maximal bounded invariant set $K_{\infty}$ in the borderline case, where we have either $\sigma_{1}=\ldots=\sigma_{l}=1$ or $\sigma_{l+1}=\ldots=\sigma_{m}=1$, remains an open problem. \hfill $\square$
\end{remark}
\begin{remark}\label{rem-res}
Let us consider the nonlinear perturbation given by $F(u) = v_0$ for $u\in X^{\alpha}$, where $v_0\in\ker A\setminus\{0\}$. We claim that the set $K_{\infty}$ is empty. Indeed, if $u$ would be a bounded full solution for the semiflow $\Phi$, then $$u(t) = S_A(t)u(0) + \int_{0}^{t} S_A(t - \tau) v_0 \,d\tau,\quad  t \ge 0,$$ which by the equality \eqref{eq-ker}, gives $$u(t) = S_A(t)u(0) + tv_0,\quad  t\ge 0.$$
Acting on this equation by the projection $Q_{0}=P_{1}+P_{2}$ and using \eqref{semig-2bb}, we obtain $$Q_{0}u(t) = S_A(t)Q_{0}u(0) + tQ_{0}v_0 = Q_{0}u(0) + tv_0, \quad t \ge 0,$$ which contradicts the assumption that the solution $u$ is bounded, because $v_0\neq 0$, and the claim follows. \hfill $\square$
\end{remark}
\begin{remark}
Let us observe that if $\sigma_{k}=0$ for some $1\le k\le m$, then the existence of the function $h_{k}\in L^{2}(\Omega)$ satisfying the corresponding inequality 
\begin{align*}
\pm f_{k}(x,u,y)\,\mathrm{sgn}\,u_{k} \ge h_{k}(x), \quad  x\in\Omega, \ u\in\R^{m}, \ y\in\R^{nm},
\end{align*}
is an obvious consequence of the assumption $(F2)$. Clearly it is enough to take $h_{k}:=\mp C_{3}$, where $C_{3}$ is the bounding constant of the maps $f_{k}$ for $1\le k\le m$. \hfill $\square$
\end{remark}
\begin{remark}
Let us assume that $I_{1},\ldots,I_{r}$ are mutually disjoint sets of natural numbers such that $I_{1}\cup\ldots\cup I_{r} = \{1,\ldots,m\}$. For any $1\le k\le r$, we define $J_{k}$ as the collection of all the indexes $i\in I_{k}$ such that the element $\sigma_{i}$ is minimal in the set $\{\sigma_{j} \ | \ j\in I_{k}\}$. Given $\epsilon_{1},\ldots,\epsilon_{r}\in \{+,-\}$, we assume that the following resonance conditions are satisfied
\begin{equation*}
\begin{split}
&(LL1)_{\epsilon_{1}} \, \left\{\begin{aligned}
&\!\sum_{k \in J_{1}}\epsilon_{1}\left(\int_{\{u_{k}>0\}}\! f_{k}^{+}(x) |u_{k}(x)|^{1-\sigma_{k}}dx \!-\!\int_{\{u_{k}<0\}} \!f_{k}^{-}(x) |u_{k}(x)|^{1-\sigma_{k}}dx \!\right)\!>\! 0 \\[3pt] & \text{for all $u_{k}\in\ker(\lambda_{k}I - A_{k})$, where $k\in J_{1}$,}\end{aligned}\right. \\[3pt]
&\hspace{100pt}\vdots\\[3pt]
&(LLr)_{\epsilon_{r}} \, \left\{\begin{aligned}
&\!\sum_{k \in J_{r}}\epsilon_{r}\left(\int_{\{u_{k}>0\}}\! f_{k}^{+}(x) |u_{k}(x)|^{1-\sigma_{k}}dx \!-\!\int_{\{u_{k}<0\}} \!f_{k}^{-}(x) |u_{k}(x)|^{1-\sigma_{k}}dx \!\right)\!>\! 0 \\[3pt] & \text{for all $u_{k}\in\ker(\lambda_{k}I - A_{k})$, where $k\in J_{r}$.}\end{aligned}\right. \\
\end{split}
\end{equation*}
Analyzing the proof of Theorems \ref{th-1} and \ref{th-2} we can check that, under the resonance conditions $(LL1)_{\epsilon_{1}}-(LLr)_{\epsilon_{r}}$ the set $K_{\infty}$ consisting of all bounded full solutions of the semiflow $\Phi$ has an admissible isolating neighborhood and the homotopy index of $K_{\infty}$ is given by
\begin{equation*}
h(\Phi, K_{\infty}) = \Sigma^{d_{\infty}(\lambda) + \beta_{1}\tilde n_{1}(\lambda)+\ldots+\beta_{r} \tilde n_{r}(\lambda)},
\end{equation*}
where, given $1\le k\le r$, we define
\begin{align*}
\tilde n_{k}(\lambda):=\sum_{j\in I_{k}}\dim\ker(\lambda_{j}I - A_{j})
\end{align*}
and furthermore we write $\beta_{k}:=1$ if $\epsilon_{k}$ has the plus sign and $\beta_{k}:=0$ otherwise. \hfill $\square$
\end{remark}
In the proof of Theorems \ref{th-1} and \ref{th-2} we use the homotopy invariance property $(H3)$ of the Rybakowski-Conley index and we deform the semiflow $\Phi$ using the following family of the differential equations
\begin{equation}\label{A-G}
\dot u(t)  = - A u(t) + H(s, u(t)), \quad  t > 0,
\end{equation}
where $H:[0,1]\times X^\alpha \to X$ is a map given by
\begin{equation}
H(s, u) := Q_{0}F(sQ_{-} u+ sQ_{+} u + Q_{0} u) + s Q_{-}F( u)+s Q_{+}F( u)
\end{equation}
for $s\in[0,1]$ and $ u\in X^\alpha$. 
\begin{remark}\label{rem-lipsch-bound}
We claim that, for any $R>0$, there is a constant $\tilde L_{R}>0$ such that 
\begin{align}\label{lok-lip-11}
\|H(s, u_{1}) - H(s, u_{2})\| \le \tilde L_{R}\| u_{1}- u_{2}\|,\quad s\in[0,1], \ \| u_{1}\|_{\alpha},\| u_{2}\|_{\alpha}\le R.
\end{align}
To check this, let us write $R':= (\|Q_{0}\|_{\alpha}+\|Q_{-}+Q_{+}\|_{\alpha}+1)R$. By the condition $(F1)$, there is a constant $L_{R'}>0$ such that 
\begin{align*}
\|F( u_{1}) - F( u_{2})\|\le L_{R'}\|  u_{1} -  u_{2}\|_{\alpha},\quad\text{if}\quad \|  u_{1}\|_{\alpha},\|  u_{2}\|_{\alpha}\le R'.
\end{align*}
Let us take $ u_{1},  u_{2}\in X^{\alpha}$ such that $\|  u_{1}\|_{\alpha},\|  u_{2}\|_{\alpha}\le R$. Then, for $k=1,2$, we have
\begin{align*}
\|sQ_{-}  u_{k}+ sQ_{+}  u_{k} + Q_{0}  u_{k}\|_{\alpha}\le (\|Q_{0}\|_{\alpha}+\|Q_{-}+Q_{+}\|_{\alpha})R = R', \ \  s\in[0,1].
\end{align*}
Therefore, for any $s\in[0,1]$, we obtain
\begin{align*}
&\|H(s,  u_{1}) - H(s,  u_{2})\| \\
& \qquad \le L_{R'}(\|Q_{0}\|\|s(Q_{-}+Q_{+})( u_{1} -  u_{2}) + Q_{0}( u_{1}- u_{2})\|_{\alpha}+\|Q_{-}+Q_{+}\|\| u_{1} -  u_{2}\|_{\alpha}) \\
&\qquad\le L_{R'}(\|Q_{0}\|+\|Q_{-}+Q_{+}\|)(\|Q_{0}\|_{\alpha}+\|Q_{-}+Q_{+}\|_{\alpha} + 1)\|  u_{1} -  u_{2}\|_{\alpha},
\end{align*}
which gives the inequality \eqref{lok-lip-11}, as desired. Observe that, by condition $(F2)$, we easily deduce that the set $\{F(v) \ | \  v\in X^{\alpha}\}$ is bounded in $X$. This implies that
\begin{equation}\label{rem-imp}
\begin{aligned}
\|H(s,  u)\| &\le \|Q_{0}\| \|F(s Q_{-}  u+s Q_{+}  u + Q_{0}  u)\| + \|Q_{-}+Q_{+}\|\|F( u)\| \\
&\le \sup\{\|F(v)\| \ | \  v\in X^{\alpha}\} ( \|Q_{0}\|+\|Q_{-}+Q_{+}\|):=C_{6},
\end{aligned}
\end{equation}
for $s\in[0,1]$ and $u\in X^\alpha$, which shows that $H$ is a bounded map. \hfill $\square$ 
\end{remark}
Arguing similarly as in Section 2, we can verify that, for any $s\in[0,1]$ and $u_{0}\in X^\alpha$, the equation \eqref{A-G} admits a unique mild solution $u(t; s, u_{0}):[0,+\infty)\to X^\alpha$ starting at $u_{0}$. In fact, it is enough to use \cite[Theorem 3.3.3]{MR610244}, \cite[Corollary 3.3.5]{MR610244} and Remark \ref{rem-lipsch-bound}. Therefore we are able to define the semiflow associated with the equation \eqref{A-G} by the following formula $$\Psi^s(t,u_{0}):=u(t;s,u_{0}),\quad t\in[0,+\infty), \ s\in[0,1], \ u_{0}\in X^\alpha.$$ 
Furthermore, applying \cite[Proposition 2.3.2]{MR1778284} once again, we deduce that the family of semiflows is continuous, that is, for any sequence $( u_{n})$ in $X^{\alpha}$ and $(s_{n})$ in $[0,1]$ such that $ u_{n}\to  u_{0}$ and $s_{n}\to s_{0}$ as $n\to\infty$, we have
\begin{align*}
\Psi^{s_{n}}(t;  u_{n})\to \Psi^{s_{0}}(t;  u_{0})\quad \text{for} \ \ t\ge 0, \ \ \text{as} \ \ n\to\infty
\end{align*}
and the convergence is uniform for the time $t$ from bounded subsets of $[0,+\infty)$. Taking into account the fact that the operator $A$ has compact resolvents, we can use \cite[Theorem 3.2.1]{MR1778284} to find that any bounded set $M\subset X^\alpha$ is {\em admissible} with respect to the family $\{\Psi^{s}\}_{s\in[0,1]}$. 

\section{Estimates for the nonlinear perturbation}
In this section we use the resonance conditions $(LL1)_{\pm}$ and $(LL2)_{\pm}$ to obtain guiding function type estimates on the nonlinear perturbation $F$. Let us observe that, due to the inclusion $L^{p}(\Omega;\R^{m})\subset L^{2}(\Omega;\R^{m})$, we can define the bilinear forms
\begin{align*}
\< u, v\>_{1} := \sum_{k=1}^{l} \int_{\Omega} u_{k}(x) v_{k}(x)\,dx, \quad\< u, v\>_{2} := \sum_{k=l+1}^{m} \int_{\Omega}  u_{k}(x) v_{k}(x)\,dx,  \  u, v\in X,
\end{align*}
that determine the functions $\| u\|_{k}^{2} := \< u, u\>_{k}$ for $u\in X$ and $k=1,2$. We intend to prove the following proposition.
\begin{proposition}\label{prop-1}
Let us assume that $\{h_{k}\}_{k=1}^{m}$ is a family of $L^{2}(\Omega)$ functions such that the inequalities $(C1)_{\pm}$ hold.
If the condition $(LL1)_{\pm}$ is satisfied, then, for any bounded set $W\subset X^{\alpha}_{-}\oplus X^{\alpha}_{+}$, there are $r>0$ and $R > 0$ such that
\begin{align}\label{ineq-f}
\pm\<F( u+ v +  w),  u\>_{1} > r 
\end{align}
for all $( u, v, w)\in N_{1}\times N_{2}\times W$ such that $\| u\|_{1}\ge R$. 
\end{proposition}
\proof Arguing by contradiction, we can suppose that there are sequences $(r_{n})$ of positive numbers, $( u_{n}, v_{n})$ in $N_1\times N_{2}$ and $(w_{n})$ in $W$ such that $r_{n}\to 0$ and $\| u_{n}\|_{1}\to\infty$ as $n\to\infty$ and furthermore 
\begin{equation}\label{c1}
\<F( u_{n} +  v_{n}+ w_{n}), u_{n}\>_{1}\le r_{n}, \quad  n \ge 1.
\end{equation}
For any $n\ge 1$, we write $ z_n :=  u_{n}/\| u_{n}\|_{1}$. Since $( z_n)$ is a bounded sequence of the finite dimensional space $N_1$ and the embedding $X^{\alpha}\subset X$ is compact, without loss of generality we can assume that there are $ z_0\in N_1$ with $\|z_0\|_{1} = 1$ and $ w_{0}\in X^{\alpha}$ such that $\| z_n -  z_0\|_{1}\to 0$ and $\| w_n -  w_0\|\to 0$ as $n\to\infty$. Let us write $ p_{n}:= u_n+ v_{n}+ w_{n}$ for $n\ge 1$. 
Then, for any $1\le k\le l$, we have $v^{k} = 0$ and 
\begin{align*}
 p^{k}_{n}/\| u_{n}\|_{1} =  z^{k}_n+  v^{k}_{n}/\| u_{n}\|_{1} +  w^{k}_{n}/\| u_{n}\|_{1} =  z^{k}_n +  w^{k}_{n}/\| u_{n}\|_{1}\to  z^{k}_{0}, \quad n\to\infty
\end{align*}
in the space $L^{2}(\Omega)$. Therefore, passing to a subsequence if necessary, we obtain 
\begin{align}\label{conv-22mm}
 w^{k}_{n}(x)\to  w^{k}_{0}(x)  \ \ \text{and} \ \ \  p^{k}_n(x)/\| u_{n}\|_{1} \to  z^{k}_{0}(x) \ \text{ as } \ n\to\infty \ \ \text{for a.a.} \ \ x\in\Omega.
\end{align}
Furthermore there are functions $a_{k},b_{k}\in L^{2}(\Omega)$, for $1\le k\le m$, such that
\begin{align}\label{est-pp-11}
| w^{k}_{n}(x)|\le a_k(x) \quad\text{and}\quad| p^{k}_n(x)/\| u_{n}\|_{1}|\le b_{k}(x),\quad x\in\Omega, \ n\ge 1.
\end{align}
By the inequality \eqref{c1}, we have the estimates
\begin{equation}
\begin{aligned}\label{gr4}
\| u_n\|^{\check\sigma_{1}-1}_{1}r_{n} &\ge \| u_n\|^{\check\sigma_{1}-1}_{1}\<F( u_{n}+ v_{n}\!+\! w_{n}), u_n\>_{1} = \| u_n\|^{\check\sigma_{1}-1}_{1}\<F( p_{n}), u_n\>_{1}\\
& = \sum_{k=1}^{l}\int_{\Omega}\| u_n\|^{\check\sigma_{1}-1}_{1}f_{k}(x, p_{n}(x),\nabla  p_{n}(x)) u^{k}_n(x) \,d x,
\end{aligned}
\end{equation}
where we recall that the number $\check\sigma_{1}$ is defined in \eqref{def-sig}. Let us assume that the inequalities $(C1)_{+}$ and resonance condition $(LL1)_{+}$ are satisfied. If we take arbitrary $1\le k\le l$, then, by \eqref{est-pp-11}, we have 
\begin{equation}
\begin{aligned}\label{est-kk-11}
&\| u_n\|^{\sigma_k-1}_{1} f_{k}(x, p_n(x),\nabla  p_{n}(x)) p^{k}_n(x)\\
&\qquad =[\mathrm{sgn}\, p_{n}^{k}(x)] f(x, p_{n}(x),\nabla  p_{n}(x))| p^{k}_{n}(x)|^{\sigma_k}| p_{n}^{k}(x)/\| u_n\|_{1}|^{1-\sigma_k}\\
& \qquad \ge h_{k}(x) | p_{n}^{k}(x)/\| u_n\|_{1}|^{1-\sigma_k}\ge -|h_{k}(x)| b_{k}(x)^{1-\sigma_k}
\end{aligned}
\end{equation}
for $x\in\Omega$ and $n\ge 1$. On the other hand, if we set $c_{0} := \sup\{\|u_n\|_{1}^{-1} \ | \ n\ge 1\}$, then, by the condition (F2) and inequality \eqref{est-pp-11}, we obtain
\begin{align*}
&-\|u_n\|^{\sigma_k-1}_{1}f_{k}(x,p_n(x), \nabla p_{n}(x))w^{k}_n(x)\ge -C_{3}c^{1-\sigma_k}_{0} a_{k}(x), \quad x\in\Omega, \ n\ge 1,
\end{align*}
which together with \eqref{est-kk-11} give the following estimate
\begin{equation}
\begin{aligned}\label{lower-major}
&\|u_n\|^{\sigma_k-1}_{1} f_{k}(x,p_n(x),\nabla p_{n}(x))u^{k}_n(x) \\
&\qquad= \|u_n\|^{\sigma_k-1}_{1} f_{k}(x,p_n(x),\nabla p_{n}(x))[p^{k}_n(x) - w^{k}_n(x)] \\
& \qquad \ge -|h_{k}(x)| b_{k}(x)^{1-\sigma_k}-C_{3}c_{0}^{1-\sigma_k}a_{k}(x)
\end{aligned}
\end{equation}
for $x\in\Omega$ and $n\ge 1$. Let us observe that in the case $\sigma_k=1$ we have 
\begin{equation*}
|h_{k}|b_{k}^{1-\sigma_k}+C_{3}c_{0}^{1-\sigma_k}a_{k} = |h_{k}|+C_{3}a_{k}\in L^{2}(\Omega)\subset L^{1}(\Omega).
\end{equation*}
Furthermore, if $\sigma_k\in[0,1)$ then $b_{k}^{1-\sigma_k}\in L^{2/(1-\sigma_k)}(\Omega)$. Since $2/(1-\sigma_k)\ge 2$ and the domain 
$\Omega\subset\R^{n}$ is bounded, it follows that $b_{k}^{1-\sigma_k}\in L^{2}(\Omega)$ and consequently
\begin{align}\label{ineq-11-aab}
\int_{\Omega} |h_{k}(x)| b_{k}(x)^{1-\sigma_k}\,dx \le \|h_{k}\|_{L^{2}(\Omega)} \|b_{k}^{1-\sigma_k}\|_{L^{2}(\Omega)} <\infty.
\end{align}
This in turn, implies that 
\begin{equation}\label{cc-55}
|h_{k}| b_{k}^{1-\sigma_k}+C_{3}c_{0}^{1-\sigma_k}a_{k}\in L^{1}(\Omega),\quad 1\le k\le l.
\end{equation}
In particular, combining \eqref{lower-major} and \eqref{cc-55} implies that
\begin{align}\label{kk-11aa}
\liminf_{n\to\infty}\|u_n\|^{\check\sigma_{1}-\sigma_{k}}_{1}\int_{\Omega}\|u_n\|^{\sigma_{k}-1}_{1}f_{k}(x,p_{n}(x),\nabla p_{n}(x))u^{k}_n(x)\,dx \ge 0
\end{align}
for all $1\le k\le m$ such that $k\not\in J_{1}$. Let us define the following sets $$I_{1}(z_{0}):=\{1\le j\le l \ | \ z_{0}^{j} = 0\},\quad I_{2}(z_{0}):=\{1\le j\le l \ | \ z_{0}^{j} \neq 0\}.$$
If we take an arbitrary $k\in I_{1}(z_{0})$, then the latter convergence of \eqref{conv-22mm} says that 
\begin{equation}\label{cc-44}
p^{k}_n(x)/\|u_{n}\|_{1}\to 0, \ \ \text{as}  \ \ n\to\infty \ \ \text{for a.a.} \ \ x\in\Omega.
\end{equation}
Observe that using the inequalities $(C1)_{+}$, we obtain the following estimates
\begin{equation}
\begin{aligned}\label{cc-11}
&\int_{\Omega}\|u_n\|^{\sigma_{k}-1}_{1}f_{k}(x,p_{n}(x),\nabla p_{n}(x))p^{k}_{n}(x)\,dx \\
&\qquad\quad =\int_{\Omega}[\mathrm{sgn}\, p_{n}^{k}(x)]f_{k}(x,p_{n}(x),\nabla p_{n}(x))|p^{k}_{n}(x)|^{\sigma_k} |p^{k}_n(x)/\|u_n\|_{1}|^{1-\sigma_k}\,dx\\
&\qquad\quad \ge\int_{\Omega}h_{k}(x)|p^{k}_n(x)/\|u_n\|_{1}|^{1-\sigma_k}\,dx,\quad k\ge 1.
\end{aligned}
\end{equation}
Furthermore, in view of \eqref{def-sig}, we have $\check\sigma_{1}<1$, which implies that at least one of the inequalities $1-\sigma_{k}>0$ and $\check\sigma_{1}-\sigma_{k}<0$ holds true. Hence, we have the estimate
\begin{align*}
|h_{k}(x)\|u_n\|^{\check\sigma_{1}-\sigma_{k}}_{1}|p^{k}_n(x)/\|u_n\|_{1}|^{1-\sigma_k}|\le  c_{0}^{\sigma_{k}-\check\sigma_{1}}|h_{k}(x)|b_{k}(x)^{1-\sigma_k},\quad x\in\Omega, \ n\ge 1,
\end{align*}
where we recall that $c_{0} := \sup\{\|u_n\|_{1}^{-1} \ | \ n\ge 1\}$ and $c_{0}^{\sigma_{k}-\check\sigma_{1}}|h_{k}|b_{k}^{1-\sigma_k}\in L^{1}(\Omega)$ as the inequality \eqref{ineq-11-aab} says.
Combining this with \eqref{cc-44}, \eqref{cc-11} and the dominated convergence theorem, we infer that
\begin{equation}
\begin{aligned}\label{cc-33}
&\liminf_{n\to\infty}\int_{\Omega}\|u_n\|^{\check\sigma_{1}-1}_{1}f_{k}(x,p_{n}(x),\nabla p_{n}(x))p^{k}_{n}(x)\,dx \\
&\qquad\quad\ge\liminf_{n\to\infty}\int_{\Omega}h_{k}(x)\|u_n\|^{\check\sigma_{1}-\sigma_{k}}_{1}|p^{k}_n(x)/\|u_n\|_{1}|^{1-\sigma_k}\,dx \\
&\qquad\quad=\lim_{n\to\infty}\int_{\Omega}h_{k}(x)\|u_n\|^{\check\sigma_{1}-\sigma_{k}}_{1}|p^{k}_n(x)/\|u_n\|_{1}|^{1-\sigma_k}\,dx = 0.
\end{aligned}
\end{equation}
Observe that using the former inequality of \eqref{est-pp-11} and the condition $(F2)$, we obtain 
\begin{align*}
\int_{\Omega}|f_{k}(x,p_n(x),\nabla p_{n}(x))w^{k}_{n}(x)|\,dx\le C_{3}\int_{\Omega}a_{k}(x)\,dx\le C_{3}|\Omega|^{1/2}\|a_{k}\|_{L^{2}(\Omega)}, \quad n\ge 1
\end{align*}
which together with the fact that $\check\sigma_{1}<1$, give
\begin{align}\label{cc-22}
\lim_{n\to\infty}\int_{\Omega}\|u_n\|^{\check\sigma_{1}-1}_{1}f_{k}(x,p_n(x),\nabla p_{n}(x))w^{k}_{n}(x)\,dx = 0.
\end{align}
Therefore, by \eqref{cc-33} and \eqref{cc-22}, for any $k\in I_{1}(z_{0})$, we obtain 
\begin{equation}
\begin{aligned}\label{eq-bb-11}
&\liminf_{n\to\infty}\int_{\Omega}\|u_n\|^{\check\sigma_{1}-1}_{1}f_{k}(x,p_{n}(x),\nabla p_{n}(x))u^{k}_{n}(x)\,dx\\
&\hspace{30pt} \ge\liminf_{n\to\infty}\int_{\Omega}\|u_n\|^{\check\sigma_{1}-1}_{1}f_{k}(x,p_{n}(x),\nabla p_{n}(x))p^{k}_{n}(x)\,dx \\
&\hspace{55pt} -\lim_{n\to\infty}\int_{\Omega}\|u_n\|^{\check\sigma_{1}-1}_{1}f_{k}(x,p_{n}(x),\nabla p_{n}(x))w^{k}_{n}(x)\,dx \ge 0.
\end{aligned}
\end{equation}
Let us take an arbitrary $k\in I_{2}(z_{0})$ and define the following sets
$$\Omega^{k}_{+}:=\{x\in\Omega \ | \ z_{0}^{k}(x) >0\}, \qquad \Omega^{k}_{-}:=\{x\in\Omega \ | \ z_{0}^{k}(x) <0\}.$$ 
If $x\in\Omega^{k}_{+}$ then $z^{k}_{n}(x)> 0$ for sufficiently large $n\ge 1$, which gives 
\begin{align}\label{conv-33-mm}
p^{k}_{n}(x) = z^{k}_n(x)\|u_n\|_{1}+w^{k}_{n}(x)\to+\infty, \ \ \text{as} \ \ n\to\infty \ \text{ for a.a. } \ x\in\Omega^{k}_{+}.
\end{align}
Observe that the similar argument shows that
\begin{align}\label{conv-33-mm22}
p^{k}_{n}(x)=z^{k}_n(x)\|u_n\|_{1}+w^{k}_{n}(x)\to-\infty, \ \ \text{as} \ \ n\to\infty \ \text{ for a.a. } \ x\in\Omega^{k}_{-}.
\end{align}
Therefore, by \eqref{conv-11-aa-kk}, \eqref{conv-22mm}, \eqref{conv-33-mm} and \eqref{conv-33-mm22}, for any $x\in\Omega^{k}_{+}$, we have
\begin{gather}\label{lim-11}
\lim_{n\to\infty}f_{k}(x,p_{n}(x),\nabla p_{n}(x))|p^{k}_{n}(x)|^{\sigma_k} |p^{k}_n(x)/\|u_n\|_{1}|^{1-\sigma_k}=f^{+}_{k}(x)|z^{k}_0(x)|^{1-\sigma_k}
\end{gather}
and furthermore, for any $x\in\Omega^{k}_{-}$, the following limit hold
\begin{gather}\label{lim-22}
\lim_{n\to\infty}f_{k}(x,p_{n}(x),\nabla p_{n}(x))|p^{k}_{n}(x)|^{\sigma_k} |p^{k}_n(x)/\|u_n\|_{1}|^{1-\sigma_k}= f^{-}_{k}(x)|z^{k}_0(x)|^{1-\sigma_k}.
\end{gather}
It is not difficult to check that for $\sigma_{k}\in(0,1]$ the limit \eqref{conv-11-aa-kk} implies that
\begin{equation*}
\lim_{s\to \pm\infty} f_{k}(x,u+se_{k},y) = 0 \quad\text{for}\quad x\in\Omega, \ u\in\mathrm{span}\,\{e_{j} \ | \ j\neq k\} \ \text{ and } \ y\in\R^{nm},
\end{equation*}
which together with \eqref{conv-22mm} and \eqref{conv-33-mm} give
\begin{align}\label{lim-33}
\lim_{n\to\infty}\|u_n\|^{\sigma_k-1}_{1}f_{k}(x,p_n(x),\nabla p_{n}(x))w^{k}_n(x)=0,\quad x\in\Omega^{k}_{+}\cup \Omega^{k}_{-}.
\end{align}
Applying \eqref{lim-11}, \eqref{lim-22} and \eqref{lim-33} to the following formula
\begin{equation}
\begin{aligned}\label{equa-pp-11bb}
&\|u_n\|^{\sigma_k-1}_{1}f_{k}(x,p_n(x),\nabla p_{n}(x))u^{k}_{n}(x) \\
&\qquad= [\mathrm{sgn}\,p_{n}^{k}(x)]f_{k}(x,p_{n}(x),\nabla p_{n}(x))|p^{k}_{n}(x)|^{\sigma_k}|p_{n}^{k}(x)/\|u_n\|_{1}|^{1-\sigma_k} \\
&\qquad\qquad - \|u_n\|^{\sigma_k-1}_{1}f_{k}(x,p_n(x),\nabla p_{n}(x))w^{k}_n(x),
\end{aligned}
\end{equation} 
we obtain the limits
\begin{align}\label{c-11}
&\lim_{n\to\infty}\|u_n\|^{\sigma_k-1}_{1}f_{k}(x,p_n(x),\nabla p_{n}(x))u^{k}_n(x)=  f^{+}_{k}(x)|z^{k}_0(x)|^{1-\sigma_k}, &&  x\in\Omega^{k}_{+},\\ \label{c-22}
&\lim_{n\to\infty}\|u_n\|^{\sigma_k-1}_{1}f_{k}(x,p_n(x),\nabla p_{n}(x))u^{k}_n(x)= - f^{-}_{k}(x)|z^{k}_0(x)|^{1-\sigma_k}, && x\in\Omega^{k}_{-}.
\end{align}
Therefore we can apply \eqref{lower-major}, \eqref{cc-55}, \eqref{c-11}, \eqref{c-22} and Fatou's lemma to obtain
\begin{align}\label{est-w-1}
\liminf_{n\to\infty}\!\int_{\Omega_{+}^{k}} \!\!\|u_n\|^{\sigma_{k}-1}_{1}\!f_{k}(x,p_n(x),\!\nabla p_{n}(x))u^{k}_n(x)dx\ge\int_{\Omega_{+}^{k}}\!\!f^{+}_{k}(x)|z^{k}_0(x)|^{1-\sigma_k}dx
\end{align}
and
\begin{align}\label{est-w-1b}
\hspace{-6pt}\liminf_{n\to\infty}\!\int_{\Omega_{-}^{k}} \!\!\|u_n\|^{\sigma_{k}\hspace{-0.8pt}-1}_{1}\!f_{k}(x,p_n(x),\!\nabla p_{n}(x))u^{k}_n(x)dx\ge\!-\!\int_{\Omega_{-}^{k}}\!\!f^{-}_{k}(x)|z^{k}_0(x)|^{1-\sigma_k}dx.
\end{align}
Observe that, for any $k\in I_{2}(z_{0})$, we have $z^{k}_{0}\in \mathrm{Ker}\,(A_{k} - \lambda_{k}I)\setminus\{0\}$ and hence,
by the unique continuation property for elliptic operators in the scalar case (see \cite[Theorem 1.1]{MR882069} and \cite[Proposition 3]{MR1151266}), we infer that the set $\Omega_{0}^{k}:=\{x\in\Omega \ | \ z_{0}^{k}(x)=0\}$ has the Lebesgue measure equal to zero. Furthermore, for any $k\not\in J_{1}$, we have the convergence $\|u_n\|^{\check\sigma_{1}-\sigma_{k}}_{1}\to 0$ as $n\to \infty$. Consequently, by \eqref{gr4}, \eqref{kk-11aa}, \eqref{eq-bb-11}, \eqref{est-w-1} and \eqref{est-w-1b}, we find that 
\begin{equation*}
\begin{aligned}
0 & \ge \sum_{k=1}^{l}\liminf_{n\to\infty}\int_{\Omega}\|u_n\|^{\check\sigma_{1}-1}_{1}f_{k}(x,p_{n}(x),\nabla p_{n}(x))u^{k}_n(x)\,dx\\
& \ge \sum_{k\in I_{2}(z_{0})}\liminf_{n\to\infty}\|u_n\|^{\check\sigma_{1}-\sigma_{k}}_{1}\int_{\Omega}\|u_n\|^{\sigma_{k}-1}_{1}f_{k}(x,p_{n}(x),\nabla p_{n}(x))u^{k}_n(x)\,dx\\
& \ge \sum_{k\in I_{2}(z_{0})\cap J_{1}}\liminf_{n\to\infty}\int_{\Omega}\|u_n\|^{\sigma_k-1}_{1}f_{k}(x,p_{n}(x),\nabla p_{n}(x))u^{k}_n(x)\,dx\\
&\ge \sum_{k\in I_{2}(z_{0})\cap J_{1}}\left(\int_{\Omega_{+}^{k}}f^{+}_{k}(x)|z^{k}_0(x)|^{1-\sigma_k}\,dx - \int_{\Omega_{-}^{k}}f^{-}_{k}(x)|z^{k}_0(x)|^{1-\sigma_k}\,dx\right)\\
&= \sum_{k\in J_{1}}\left(\int_{\Omega_{+}^{k}}f^{+}_{k}(x)|z^{k}_0(x)|^{1-\sigma_k}\,dx - \int_{\Omega_{-}^{k}}f^{-}_{k}(x)|z^{k}_0(x)|^{1-\sigma_k}\,dx\right),
\end{aligned}
\end{equation*}
where the last equality follows from the fact that $1>\check\sigma_{1} = \sigma_k$ for $k\in J_{1}$ and from the definition of the set $I_{1}(z_{0})$. This in turn contradicts $(LL1)_{+}$ and proves the plus sign case of the inequality \eqref{ineq-f}. Analogous argument, the details of which we leave to the reader, shows that the inequalities $(C1)_{-}$ and resonance condition $(LL1)_{-}$ imply the minus sign form of \eqref{ineq-f} and thus the proof of the proposition is completed. \hfill $\square$\\

Following the lines of the above proof, we can show the following proposition concerning the guiding function type estimates for the nonlinear perturbation $F$, with respect to the space $N_{2}$.
\begin{proposition}\label{prop-2}
Suppose that $\{h_{k}\}_{k=1}^{m}$ is a family of $L^{2}(\Omega)$ functions such that the inequalities $(C2)_{\pm}$ hold. If the condition $(LL2)_{\pm}$ is satisfied, then, for any bounded set $W\subset X^{\alpha}_{-}\oplus X^{\alpha}_{+}$, there is $r>0$ and $R > 0$ such that
\begin{align*}
\pm\<F(u+v + w), v\>_{2} > r 
\end{align*}
for all $(u, v, w)\in N_{1}\times N_{2}\times W$ such that $\|v\|_{2}\ge R$. 
\end{proposition}

\section{Estimates for bounded solutions of the homotopy flow}
 We begin with the following lemma, which provides some {\em a priori} bounds for the projections of solutions of the parametrized equation \eqref{A-G} onto the space $X_{-}^{\alpha}\oplus X_{+}^{\alpha}$. 
\begin{proposition}\label{prop:3}
There is $R_{0} > 0$ such that for any $s\in[0,1]$ and for any bounded full solution $u=u_s$ for the semiflow $\Psi^s$, the following inequality holds
\begin{equation}\label{ineq-ppjjhh}
\|Q_{k}u(t)\|_\alpha\le R_{0} \quad\text{for} \ \  t\in\R \ \ \text{and} \ \ k\in\{-,+\}.
\end{equation}
\end{proposition}
\noindent\textbf{Proof.} Let us assume that $u=u_s$ is a full solution for the equation \eqref{A-G} for some $s\in[0,1]$. From the continuity of the projections $Q_{-},Q_{+}:X^\alpha\to X^\alpha$ we deduce the boundedness of the sets $\{Q_{-}u(t) \ | \ t\ge 0\}$ and $\{Q_{+}u(t) \ | \ t\le 0\}$ in the space $X^\alpha$. Since $u$ is a solutions of the semiflow $\Psi^{s}$, it follows that the equality $\Psi^s(t-t',u(t')) = u(t)$ holds for all $t,t'\in\R$ such that $t\ge t'$.  This in turn can be written in the following integral form
\begin{equation}\label{ful-sol}
u(t) =  S_A(t - t')u(t') + \int_{t'}^{t} S_A(t - \tau)H(s,u(\tau)) \,d \tau, \ \ t>t'.
\end{equation}
Acting on \eqref{ful-sol} by the operator $Q_{k}$, where $k\in\{-,+\}$, and using \eqref{semig-2}, we obtain
\begin{equation}\label{fs2}
Q_{k}u(t) = S_A(t - t')Q_{k}u(t') + \int_{t'}^t S_A(t - \tau)Q_{k}H(s,u(\tau)) \,d\tau, \ \ t>t'.
\end{equation}
Since the semigroup $\{S_A(t)\}_{t\ge 0}$ extends on the space $X_{-}$ to the $C_0$ group of bounded linear operators, we can apply $S_{A}(t' - t)$ on the formula \eqref{fs2} to derive
\begin{equation}\label{int-form-1}
S_A(t' - t) Q_{-}u(t) = Q_{-}u(t') + \int_{t'}^t S_A(t' - \tau)Q_{-}H(s,u(\tau)) \d \tau, \ \ t\ge t'.
\end{equation}
Then, the inequalities \eqref{cont-emb} and \eqref{bb-nn-2} imply that
\begin{align*}
\|S_A(t' - t) Q_{-}u(t)\| \le C_{5} e^{c(t' - t)} \|Q_{-}u(t)\|\le C_{2}C_{5} e^{c(t' - t)} \|Q_{-}u(t)\|_{\alpha}
\end{align*}
and hence, using the boundedness of the set $\{Q_{-}u(t) \ | \ t\ge 0\}$ in $X^{\alpha}$, we find that
\begin{equation}\label{limit-11}
\|S_A(t' - t) Q_{-}u(t)\| \to 0,\quad t\to +\infty.
\end{equation}
Combining \eqref{int-form-1} with \eqref{bb-nn-2} and \eqref{rem-imp}, we obtain
\begin{equation}
\begin{aligned}\label{ineq-33-55}
\|Q_{-}u(t')\| & \le \|S_A(t' - t) Q_{-}u(t)\| + C_{5}\int_{t'}^t e^{c (t' - \tau)}\|Q_{-}H(s,u(\tau))\| \,d \tau \\
& \le \|S_A(t' - t) Q_{-}u(t)\| + C_{5}C_{6} \|Q_{-}\| \int_{t'}^t e^{c (t' - \tau)} \,d \tau \\
& \le \|S_A(t' - t) Q_{-}u(t)\| + C_{5}C_{6} c^{-1}\|Q_{-}\|.
\end{aligned}
\end{equation}
Since $X_{-}$ is a finite dimensional space there is constant $C_{7} > 0$ such that 
\begin{equation}\label{est-c-3}
\| u\|_\alpha \le C_{7}\| u\|,\quad  u\in X_{-}.
\end{equation}
Passing in \eqref{ineq-33-55} to the limit with $t\to +\infty$ and using \eqref{limit-11} with \eqref{est-c-3} we obtain
\begin{align}\label{fs7}
\|Q_{-}u(t')\|_{\alpha} \le C_{5}C_{6}C_{7}\|Q_{-}\|c^{-1},\quad t'\in\R,
\end{align}
which gives the desired estimate \eqref{ineq-ppjjhh} for $k=1$. Let us observe that, by the inequalities \eqref{cont-emb} and \eqref{bb-nn-1}, we have
\begin{align*}
 \|S_A(t - t')Q_{+}u(t')\|_\alpha & \le C_{5}e^{- c (t - t')}(t - t')^{-\alpha} \, \|Q_{+}u(t')\| \\
&\le  C_{2}C_{5}e^{- c (t - t')}(t - t')^{-\alpha} \, \|Q_{+}u(t')\|_\alpha,
\end{align*}
which together with the boundedness of the set $\{Q_{+}u(t) \ | \ t\le 0\}$ in $X^{\alpha}$ give
\begin{equation}\label{fs3}
\|S_A(t - t')Q_{+}u(t')\|_\alpha \to 0,\quad t'\to - \infty.
\end{equation}
Combining the formula \eqref{fs2} with the inequalities \eqref{bb-nn-1} and \eqref{rem-imp}, we obtain
\begin{equation}\label{row-now2}
\begin{aligned}
\|Q_{+}u(t)\|_{\alpha} &\le \|S_A(t - t')Q_{+}u(t')\|_\alpha + C_{5}\int_{t'}^{t}\frac{e^{-c(t-\tau)}}{(t - \tau)^\alpha}\,\|Q_{+}H(s,u(\tau))\|\,d\tau\\
& \le \|S_A(t - t')Q_{+}u(t')\|_\alpha + C_{5}C_{6} \|Q_{+}\| \int_{t'}^t  \frac{e^{- c (t - \tau)}}{(t - \tau)^\alpha} \,d\tau.
\end{aligned}
\end{equation}
If we take $t,t'\in\R$ with $t'+1<t$, then we have the following estimates
\begin{equation*}
\begin{aligned}
& \int_{t'}^{t} \frac{e^{- c (t - \tau)}}{(t - \tau)^{\alpha}}\,d\tau = \int_{t'}^{t-1} \frac{e^{- c (t - \tau)}}{(t - \tau)^{\alpha}}\,d\tau + 
\int_{t-1}^{t} \frac{e^{- c (t - \tau)}}{(t - \tau)^{\alpha}}\,d\tau \\
&\qquad\le\int_{t'}^{t-1} e^{- c (t - \tau)}\,d\tau + \int_{t-1}^{t} (t - \tau)^{-\alpha}\,d\tau\le e^{-c}/c + 1/(1 - \alpha)
\end{aligned}
\end{equation*}
that together with the inequality \eqref{row-now2}, provide
\begin{equation}\label{equa-lim-1}
\|Q_{+}u(t)\|_\alpha \le \|S_A(t - t')Q_{+}u(t')\|_\alpha + C_{5}C_{6}\|Q_{+}\|(e^{-c}/c + 1/(1 - \alpha)).
\end{equation}
Using \eqref{fs3} and passing in \eqref{equa-lim-1} to limit with $t'\to -\infty$ we infer that 
\begin{equation}\label{fs9}
\|Q_{+}u(t)\|_\alpha \le C_{5}C_{6}\|Q_{+}\|\left(e^{-c}/c  + 1/(1 - \alpha)\right),\quad t\in\R.
\end{equation}
Thus we obtain the estimate \eqref{ineq-ppjjhh} and the proof of the proposition is completed. \hfill $\square$ \\

We proceed to the following proposition, which provides the estimates for the projections of the solutions of the equation \eqref{A-G} onto the space $N_{1}$. 

\begin{proposition}\label{prop:3a}
Let us assume that $W\subset X^{\alpha}_{-}\oplus X^{\alpha}_{+}$ is a ball centered at the origin and $r, R > 0$ are such that either
\begin{align}\label{ineq-kk-hh-11}
\hspace{-15pt}\<F(u+v + w), u\>_{1} > r \ \text{ for } \ (u,v,w)\in N_{1}\times N_{2}\times W \ \text{ with } \ \|u\|_{1}\ge R
\end{align}
or
\begin{align}\label{ineq-kk-hh-22}
\hspace{-5pt}\<F(u+v + w), u\>_{1} < -r \ \text{ for } \ (u,v,w)\in N_{1}\times N_{2}\times W \ \text{ with } \ \|u\|_{1}\ge R. 
\end{align}
Then, for any $s\in[0,1]$ and any bounded full solution $u$ of the semiflow $\Psi^s$ such that $Q_{-}u(t)+Q_{+}u(t)\in W$ for $t\in\R$, the following inequality holds
\begin{equation}\label{b-sol-2}
\|P_{1}u(t)\|_{1}\le R,\quad t\in\R.
\end{equation}
\end{proposition}
\proof Since $N_{1}$ is finite dimensional space, the functions $\|\hspace{0.7pt}\cdot\hspace{0.7pt}\|_{1}$ and $\|\hspace{0.7pt}\cdot\hspace{0.7pt}\|$ are equivalent norms on $N_{1}$. Hence the boundedness of the solution $u$ in the space $X$ gives
\begin{align}\label{b-sol}
\sup_{t\in\R} \|P_{1}u(t)\|_{1}<+\infty.
\end{align}
We argue by a contradiction and assume that there is $s\in[0,1]$ and a full solution $u$ of the semiflow $\Psi^{s}$ such that $\|P_{1}u(t_{0})\|_{1}>R$ for some $t_{0}\in\R$. Acting by the operator $P_{1}$ on the integral formula
\begin{equation*}
u(t) = S_A(t - t')u(t') + \int_{t'}^{t} S_A(t - \tau)H(s,u(\tau))\,d\tau, \quad t\ge t'
\end{equation*}
and using \eqref{semig-2bb} we obtain
\begin{equation*}
P_{1}u(t) = S_A(t - t')P_{1}u(t') + \int_{t'}^{t} S_A(t - \tau)P_{1}H(s,u(\tau))\,d\tau,
\end{equation*}
which by the kernel equality \eqref{eq-ker}, takes the following form
\begin{equation*}
P_{1}u(t) = P_{1}u(t') + \int_{t'}^t P_{1}F(sQ_{-}u(\tau)+sQ_{+}u(\tau) + P_{1}u(\tau)+ P_{2}u(\tau))\,d\tau.
\end{equation*}
Let us assume that $\<\cdot,\cdot\>$ is the standard scalar product on $L^{2}(\Omega,\R^{m})$. In view of the fact that the spaces $N_{1}$, $N_{2}$, $X_{-}$ and $X_{+}$ are mutually orthogonal, for any $u\in N_{1}$, $v\in N_{2}$ and $w\in X^{\alpha}_{-}\oplus X^{\alpha}_{+}$, we have
\begin{equation*}
\!\begin{aligned}
&\<P_{1}F(u+v + w), u\>_{1} = \<P_{1}F(u+v + w), u\> = \<(P_{1}+P_{2})F(u+v + w), u\> \\
&\qquad = \<Q_{0}F(u+v + w), u\> = \<F(u+v + w), u\> = \<F(u+v + w), u\>_{1},
\end{aligned}
\end{equation*}
which in turn, for any $t\in\R$, yields 
\begin{equation}
\begin{aligned}\label{eq-ff-11}
&\frac{d}{dt} \|P_{1}u(t)\|_{1}^{2} = 2\< \frac{d}{dt}P_{1} u(t), P_{1}u(t)\>_{1} \\
&\quad = 2\<P_{1}F(sQ_{-}u(t)+sQ_{+}u(t) + P_{1}u(t)+P_{2}u(t)), P_{1}u(t)\>_{1}\\
&\quad = 2\<F(sQ_{-}u(t)+ sQ_{+}u(t) + P_{1}u(t)+P_{2}u(t)), P_{1}u(t)\>_{1}.
\end{aligned}
\end{equation}
Let us assume that the condition \eqref{ineq-kk-hh-11} is satisfied. If we define 
$$t^{+}_{0}:=\sup\,\{t\ge t_{0} \ | \ \|u(t')\|_{1}\ge R \,\text{ for }\, t'\in[t_{0},t]\},$$ then $t_{0}^{+}=+\infty$, because otherwise, the fact that $\|u(t_{0}^{+})\|_{1}\ge R$ together with \eqref{eq-ff-11} and the inequality \eqref{ineq-kk-hh-11} would yield
\begin{align*}
\frac{d}{dt} \|P_{1}u(t)\|_{1}^{2}\, |_{t=t^{+}_{0}} = 2\<F(s(Q_{-}+Q_{+})u(t^{+}_{0}) + P_{1}u(t^{+}_{0})+P_{2}u(t^{+}_{0})), P_{1}u(t^{+}_{0})\>_{1}> 2r.
\end{align*}
Consequently, there would exists $\delta>0$ such that 
\begin{align*}
\|P_{1}u(t)\|_{1}> \|P_{1}u(t_{0}^{+})\|_{1}\ge R,\quad t\in[t^{+}_{0},t^{+}_{0}+\delta],
\end{align*}
contrary to the definition of the number $t^{+}_{0}$. This implies that $\|P_{1}u(t)\|_{1} \ge R$ for $t\ge t_{0}$ and hence, using \eqref{ineq-kk-hh-11} once again, we obtain
\begin{align*}
\frac{d}{dt} \|P_{1}u(t)\|_{1}^{2} = 2\<F(s(Q_{-}+Q_{+})u(t) + P_{1}u(t)+P_{2}u(t)), P_{1}u(t)\>_{1}> 2r,\quad t\ge t_{0}.
\end{align*}
Therefore the following inequality is satisfied 
\begin{align*}
\|P_{1}u(t)\|_{1}^{2}\ge \|P_{1}u(t_{0})\|_{1}^{2} + 2(t-t_{0})r, \quad t\ge t_{0},
\end{align*}
which contradicts \eqref{b-sol} and proves the estimate \eqref{b-sol-2}. On the other hand, if the condition \eqref{ineq-kk-hh-22} holds, then we write
\begin{align*}
t^{-}_{0}:=\inf\,\{t\le t_{0} \ | \ \|u(t')\|_{1}\ge R \ \text{ for } \ t'\in[t,t_{0}]\},
\end{align*}
If $t_{0}^{-}$ would be a finite real number, then the inequality $\|u(t_{0}^{-})\|_{1}\ge R$ combined with \eqref{eq-ff-11} and \eqref{ineq-kk-hh-22} would give
\begin{align*}
\frac{d}{dt} \|P_{1}u(t)\|_{1}^{2}\,|_{t=t^{-}_{0}} = 2\<F(s(Q_{-}+Q_{+})u(t^{-}_{0}) + P_{1}u(t^{-}_{0})+P_{2}u(t^{-}_{0})), P_{1}u(t^{-}_{0})\>_{1}< -2r.
\end{align*}
This in turn would imply the existence of $\delta>0$ such that 
\begin{align*}
\|P_{1}u(t)\|_{1}> \|P_{1}u(t_{0}^{-})\|_{1}\ge R,\quad t\in[t^{-}_{0}-\delta,t^{-}_{0}],
\end{align*}
which is impossible due to the definition of $t^{-}_{0}$. Therefore $t_{0}^{-}=-\infty$ and consequently $\|P_{1}u(t)\|_{1} \ge R$ for $t\le t_{0}$. Combining this inequality with \eqref{ineq-kk-hh-22} we obtain
\begin{align*}
\frac{d}{dt} \|P_{1}u(t)\|_{1}^{2} = 2\<F(s(Q_{-}+Q_{+})u(t) + P_{1}u(t)+P_{2}u(t)), P_{1}u(t)\>_{1}< -2r,\quad t\le t_{0},
\end{align*}
which after integration gives
\begin{align*}
\|P_{1}u(t_{0})\|_{1}^{2}+2(t_{0} - t)r < \|P_{1}u(t)\|_{1}^{2}, \quad t\le t_{0}.
\end{align*}
This again contradicts \eqref{b-sol} and shows the estimate \eqref{b-sol-2}. Thus the proof of the proposition is completed. \hfill $\square$ \\

In the similar way, we can show the following proposition concerning the estimates of the solutions of the equation \eqref{A-G} after projection onto the space $N_{2}$. 

\begin{proposition}\label{prop-hh-pp}
Let us assume that $W\subset X^{\alpha}_{-}\oplus X^{\alpha}_{+}$ is a ball centered at the origin and $r,R > 0$ are such that either
\begin{align*}
\hspace{-15pt}\<F(u+v + w), u\>_{2} > r \ \text{ for } \ (u,v,w)\in N_{1}\times N_{2}\times W \ \text{ with } \ \|u\|_{2}\ge R
\end{align*}
or
\begin{align*}
\hspace{-5pt}\<F(u+v + w), u\>_{2} < -r \ \text{ for } \ (u,v,w)\in N_{1}\times N_{2}\times W \ \text{ with } \ \|u\|_{2}\ge R. 
\end{align*}
Then, for any $s\in[0,1]$ and any bounded full solution $u$ of the semiflow $\Psi^s$ such that $Q_{-}u(t)+Q_{+}u(t)\in W$ for $t\in\R$, the following inequality holds
\begin{equation*}
\|P_{2}u(t)\|_{2}\le R,\quad t\in\R.
\end{equation*}
\end{proposition}

\section{Proof of Theorems \ref{th-1} and \ref{th-2}.} \noindent\textbf{Step 1.} Proposition \ref{prop:3} says that there is a constant $R_0 > 0$ such that, for any $s\in[0,1]$ and for any bounded full solution $u$ of the semiflow $\Psi^s$, we have
\begin{equation}\label{rrow}
\|Q_{-}u(t)+Q_{+}u(t)\|_\alpha\le R_0, \quad t\in\R.
\end{equation}
Let us define $W := \{u\in X^{\alpha}_{-}\oplus X^{\alpha}_{+} \ | \ \|u\|_\alpha \le R_0\}$. If condition $(LL1)_{\pm}$ is satisfied, then Proposition \ref{prop-1} asserts the existence of $r_{1}>0$ and $R_{1} > 0$ such that
\begin{equation}\label{rrow2}
\pm\<F(u+v+w), u\>_{1} > r_{1} \ \text{ for } \  (u,v,w)\in N_{1}\times N_{2}\times B \ \text{ with } \ \|u\|_1\ge R_{1}.
\end{equation}
Consequently, by the inequality \eqref{rrow} and Proposition \ref{prop:3a}, we infer that
\begin{equation}\label{est-gg-tt-11}
\|P_{1}u(t)\|_{1}\le R_{1}, \quad t\in\R.
\end{equation}
On the other hand, if condition $(LL2)_{\pm}$ hold, then Proposition \ref{prop-2} says that, there are $r_{2}>0$ and $R_{2}>0$ such that
\begin{equation}\label{rrow3}
\pm\<F(u+v+w), v\>_{2} > r_{2} \ \text{ for } \  (u,v,w)\in N_{1}\times N_{2}\times B \ \text{ with } \ \|u\|_{2}\ge R_{2}.
\end{equation}
Therefore, using the inequality \eqref{rrow} and Proposition \ref{prop-hh-pp}, we deduce that 
\begin{equation}\label{est-gg-tt-22}
\|P_{2}u(t)\|_{2}\le R_{2}, \quad t\in\R.
\end{equation}
Let us define the following sets
\begin{gather*}
M_{0} := \{u\in X^{\alpha}_{-}\oplus X^{\alpha}_{+} \ | \ \|u\|_\alpha \le R_0+1\}, \\
M_1 := \{u\in N_{1} \ | \ \|u\|_1 \le R_{1}+1\}, \quad M_2 := \{v\in N_{2} \ | \ \|v\|_2 \le R_{2}+1\}
\end{gather*}
and write $M:= M_{0}\oplus M_{1}\oplus M_2$. By the estimates \eqref{rrow}, \eqref{est-gg-tt-11} and \eqref{est-gg-tt-22}, we deduce that, if $u$ is a bounded full solution of $\Psi^{s}$, where $s\in[0,1]$, then it is contained in the interior of the set $M$, which in particular, is an admissible isolating neighborhood for the family of the semiflows $\{\Psi^s\}_{s\in[0,1]}$ and $\inv M = K_{\infty}$. Hence, by the homotopy invariance of the Rybakowski-Conley index (see property $(H3)$), we obtain
\begin{equation}\label{h1}
h(\Phi, K_{\infty}) = h(\Psi^1, K_1) = h(\Psi^0, K_0),
\end{equation}
where we denote $K_s:=\inv (M, \Psi^{s})$ for $s\in[0,1]$.\\[3pt]
\noindent\textbf{Step 2.} Let $\psi_{1}:[0,+\infty)\times X_{0}\to X_{0}$ be the semiflow associated with the equation
\begin{equation*}
\dot u(t)  = Q_{0}F(u(t)), \quad  t > 0.
\end{equation*}
We show that $M_1\oplus M_{2}$ is an isolating block for $\psi_{1}$ and its exit set is such that 
\begin{equation}\label{przypadkiblok}
(M_1\oplus M_{2})^{-}:=\!\left\{\begin{aligned}
&\partial_{X_{0}}(M_1\oplus M_{2})&& \!\!\text{ if $(C1)_{+}$, $(C2)_{+}$, $(LL1)_{+}$, $(LL2)_{+}$ hold}, \\
&(\partial_{N_{1}} M_1)\oplus M_{2}&& \!\!\text{ if $(C1)_{+}$, $(C2)_{-}$, $(LL1)_{+}$, $(LL2)_{-}$ hold}, \\
& M_1\oplus (\partial_{N_{2}}M_{2}) && \!\!\text{ if $(C1)_{-}$, $(C2)_{+}$, $(LL1)_{-}$, $(LL2)_{+}$ hold}, \\
& \emptyset && \!\!\text{ if $(C1)_{-}$, $(C2)_{-}$, $(LL1)_{-}$, $(LL2)_{-}$ hold}.
\end{aligned}\right.
\end{equation}
To this end, let us assume that $u:[-\delta_2,\delta_1) \to X_0$, where $\delta_1 > 0$, $\delta_2 \ge 0$, is a solution for $\psi_{1}$ such that $u_{0}:=u(0)\in\partial_{X_{0}}(M_1\oplus M_{2})$. Then we have
$$u(t) = u(0) + \int_0^t Q_{0}F(u(\tau)) \,d\tau, \quad t\in [-\delta_2,\delta_1),$$
which implies that, for any $t\in (-\delta_{2},\delta_{1}]$ and $k\in\{1,2\}$, we have
\begin{equation}
\begin{aligned}\label{cases-1}
\frac{1}{2}\frac{d}{dt}\|P_{k}u(t)\|_{k}^2 & = \<\frac{d}{dt} P_{k}u(t), P_{k}u(t)\>_k = \<P_{k}Q_{0}F(u(t)), P_{k}u(t)\>_{k} \\
&=\<Q_{0}F(u(t)), P_{k}u(t)\>_{k} = \<Q_{0}F(u(t)), P_{k}u(t)\> \\ 
&= \<F(u(t)), P_{k}u(t)\> =\<F(P_{1}u(t) + P_{2}u(t)), P_{k}u(t)\>_k.
\end{aligned}
\end{equation}
Let us assume that the conditions $(LL1)_{\pm}$ and $(LL2)_{\pm}$ are satisfied. Since $u_{0}\in \partial_{X_{0}}(M_1\oplus M_{2})$, we have either $u_{0}\in (\partial_{N_{1}}M_1)\oplus M_{2}$ or $u_{0}\in M_1\oplus (\partial_{N_{2}}M_{2})$. In the former case we have $\|P_{1}u_{0}\|_{1} = R_{1}+1$ and $\|P_{2}u_{0}\|_{2} \le R_{2}+1$, which together with the equation \eqref{cases-1} and inequality \eqref{rrow2}, implies that
\begin{equation}\label{ineq-e-11}
\pm\frac{d}{dt}\|P_{1}u(t)\|_1^2\,|_{t=0} > 0.
\end{equation}
On the other hand, if $u_{0}\in M_1\oplus (\partial_{N_{2}}M_{2})$, then $\|P_{1}u_{0}\|_{1} \le R_{1}+1$ and $\|P_{2}u_{0}\|_{2} = R_{2}+1$. Applying the equation \eqref{cases-1} together with the inequality \eqref{rrow3}, we obtain
\begin{equation}\label{ineq-e-22}
\pm\frac{d}{dt}\|P_{2}u(t)\|_2^2\,|_{t=0} > 0.
\end{equation}
Combining \eqref{ineq-e-11} and \eqref{ineq-e-22} we infer that the exit set of $M_1\oplus M_{2}$ takes the form
\begin{equation}\label{przypadkiblok-1}
(M_1\oplus M_{2})^{-}:=\left\{\begin{aligned}
&\partial_{X_{0}}(M_1\oplus M_{2})&& \!\text{if $(C1)_{+}$, $(C2)_{+}$, $(LL1)_{+}$, $(LL2)_{+}$ hold}, \\
& \emptyset && \!\text{if $(C1)_{-}$, $(C2)_{-}$, $(LL1)_{-}$, $(LL2)_{-}$ hold}.
\end{aligned}\right.
\end{equation}
On the other hand, if we assume that the conditions $(LL1)_{\pm}$ and $(LL2)_{\mp}$ are satisfied, then, proceeding in the same way, we obtain
\begin{equation*}
\pm\frac{d}{dt}\|P_{1}u(t)\|_1^2\,|_{t=0} > 0, \quad \text{if} \ \ u_{0}\in (\partial_{N_{1}}M_1)\oplus M_{2},
\end{equation*}
and furthermore
\begin{equation*}
\mp\frac{d}{dt}\|P_{2}u(t)\|_2^2\,|_{t=0} > 0, \quad\text{if} \ \  u_{0}\in M_1\oplus (\partial_{N_{2}}M_{2}).
\end{equation*}
This in turn implies that
\begin{equation*}
(M_1\oplus M_{2})^{-}:=\left\{\begin{aligned}
&(\partial_{N_{1}} M_{1})\oplus M_{2}&& \!\text{if $(C1)_{+}$, $(C2)_{-}$, $(LL1)_{+}$, $(LL2)_{-}$ hold}, \\
& M_1\oplus (\partial_{N_{2}} M_{2}) && \!\text{if $(C1)_{-}$, $(C2)_{+}$, $(LL1)_{-}$, $(LL2)_{+}$ hold}.
\end{aligned}\right.
\end{equation*}
which together with \eqref{przypadkiblok-1} provide \eqref{przypadkiblok} as desired. In particular, we infer that $M_{1}\oplus M_{2}$ is an isolated neighborhood for the invariant set 
$K_0^1:=\inv(\psi_1, M_{1}\oplus M_{2})$ and its homotopy index is given by 
\begin{equation}\label{eq-hh-dd-11}
h(\psi_1, K_0^1) = \left\{\begin{aligned}
&\Sigma^{n_{1}(\lambda)+n_{2}(\lambda)}&& \text{if $(C1)_{+}$, $(C2)_{+}$, $(LL1)_{+}$, $(LL2)_{+}$ hold},\\
&\Sigma^{n_{1}(\lambda)}&& \text{if $(C1)_{+}$, $(C2)_{-}$, $(LL1)_{+}$, $(LL2)_{-}$ hold},\\
&\Sigma^{n_{2}(\lambda)}&& \text{if $(C1)_{-}$, $(C2)_{+}$, $(LL1)_{-}$, $(LL2)_{+}$ hold},\\
&\Sigma^{0}&& \text{if $(C1)_{-}$, $(C2)_{-}$, $(LL1)_{-}$, $(LL2)_{-}$ hold}.
\end{aligned}\right.
\end{equation}
\noindent\textbf{Step 3.} Let us assume that $\psi_2$ is a semiflow obtained by the restriction of the semigroup $\{S_{A}(t)\}_{t\ge 0}$ to the space $X^{\alpha}_{-}\oplus X^{\alpha}_{+}$, that is, $$\psi_2(t, u):= S_A(t)  u,\quad t\ge 0, \  u\in X^{\alpha}_{-}\oplus X^{\alpha}_{+}.$$ Combining the estimates \eqref{bb-nn-2} and \eqref{bb-nn-3} with the following commutative property $$(\delta I  + A)^\alpha S_A(t) u = S_A(t) (\delta I  + A)^\alpha  u, \quad  u\in X^\alpha,$$ we deduce that 
\begin{equation*}
\begin{aligned}
\|S_A(t) u\|_\alpha &\le C_{5} e^{c t}\| u\|_\alpha, && t\le 0, \  u\in X_{-},\\
\|S_A(t) u\|_\alpha &\le C_{5} e^{- c t}\| u\|_\alpha, && t\ge 0, \  u\in X_{+}.
\end{aligned}
\end{equation*}
Hence \cite[Theorem 11.1]{MR910097} shows that $M_0$ is an admissible isolating neighborhood, $K_0^2:=\inv(\psi_2, M_{0})=\{0\}$ and the homotopy index of $K_{0}^{2}$ satisfies 
\begin{equation}\label{con-deg}
    h(\psi_2, K_0^2) = \Sigma^{\dim X_{-}} = \Sigma^{d_{\infty}(\lambda)}.
\end{equation}
where the last equality is a consequence of Remark \ref{rem-dim-eq}. \\[3pt]
\noindent\textbf{Step 4.} Let us observe that the semiflow $\Psi^0$ corresponding to the equation
\begin{equation*}
\dot u(t)  = - A u(t) + Q_{0}F(Q_{0}u(t)), \quad  t > 0,
\end{equation*}
satisfies the following equality
\begin{equation*}
\Psi^0(t,u + v)  = \psi_1(t,u)+\psi_2(t,v),\quad  u\in X_0, \  v\in X^{\alpha}_{-}\oplus X^{\alpha}_{+}.
\end{equation*}
Therefore, by the multiplication property of the homotopy index $(H2)$, we have
\begin{equation}\label{mult-11}
h(\Psi^0, K_0) = h(\psi_1, K_0^1)\wedge h(\psi_2, K_0^2).
\end{equation}
Combining \eqref{h1}, \eqref{mult-11} and \eqref{con-deg} we deduce that 
\begin{align*}
h(\Phi, K_{\infty}) = h(\psi_1, K_{0}^{1})\wedge\Sigma^{d_{\infty}(\lambda)},
\end{align*}
which together with \eqref{eq-hh-dd-11} provides the homotopy index formulas \eqref{form-conley-1}, \eqref{form-conley-1b}, \eqref{form-conley-2} and \eqref{form-conley-2b}. Thus the proof of Theorems \ref{th-1} and \ref{th-2} is completed. \hfill $\square$

\section{Applications to the existence of connecting solutions}
In this section we provide applications of Theorems \ref{th-1} and \ref{th-2} to study the existence of solutions connecting stationary points for the system \eqref{sys-res} with the assumption that the operator $A$ has the following particular form 
\begin{align*}
Au=(A_{0}u_{1},\ldots,A_{0}u_{m}) + (\lambda_{1}u_{1},\dots,\lambda_{m}u_{m}), \ u=(u_{1},\ldots,u_{m})\in C^{2}(\o\Omega;\R^{m}),
\end{align*}
where $\lambda_{1},\lambda_{2},\ldots,\lambda_{m}$ are real numbers and $A_{0}$ is the Laplace operator with the Dirichlet boundary conditions, that is
\begin{equation*}
\left\{\begin{aligned}
D(A_{0})&:= \mathrm{cl}_{W^{2,p}(\Omega)}\{u_{0}\in C^{2}(\o\Omega) \ | \ u_{0}(x) = 0 \ \text{ for } \ x\in\partial\Omega\},\\
A_{0}u_{0}&:= -\Delta u_{0},\quad u_{0}\in D(A_{0}).
\end{aligned}\right.
\end{equation*}
Here we recall that $p\ge 2n$ is the real number fixed in \eqref{alpha-p}.
Furthermore, we require that, for any $1\le k\le m$, the map $f_{k}:\o\Omega\times\R^{m}\times\R^{nm}\to\mathbb{R}$ is continuously differentiable and has the following property 
\begin{equation}\label{f-cond}
f(x,0,0) = 0,\quad D_{u}f(x,0,0) = G, \quad D_{y}f(x,0,0) = 0, \quad x\in\Omega,
\end{equation}
where $G$ is a symmetric $m\times m$ matrix, which is not dependent on $x\in\Omega$. 
Then, it is not difficult to check that the regularity assumption for $f$ implies that $F:X^{\alpha}\to X$ is a $C^{1}$ map satisfying conditions $(F1)$ and $(F2)$, where we recall that $\alpha\in(3/4,1)$ as it is fixed in \eqref{alpha-p}. Furthermore, from the assumption \eqref{f-cond} we infer that $F(0) = 0$ and the derivative of $F$ at the origin is a linear map $DF(0):X^{\alpha}\to X$ that, for any $w=(w_{1},w_{2},\ldots,w_{m})\in X^{\alpha}$, is given by 
\begin{equation}\label{der-11}
DF(0)[w](x) = Gw(x)\quad \text{for a.a.} \ \ x\in\Omega.
\end{equation}
If we define $\Lambda:=\mathrm{diag}\,(\lambda_{1},\lambda_{2},\ldots,\lambda_{m})$, then the matrix $G+\Lambda$ is symmetric and consequently its spectrum $\sigma(G+\Lambda)$ consists of a finite sequence of real eigenvalues $\theta_{1}\le \theta_{2}\le\ldots\le \theta_{m}$. Let us define 
\begin{align*}
d_{0}(\lambda) := \sum_{k=1}^{m}\sum_{\nu<\theta_{k}} \dim\ker(\nu I- A_{0}),
\end{align*}
where we write $\lambda:=(\lambda_{1},\lambda_{2},\ldots,\lambda_{m})$ and, in the above summation, the parameter $\nu$ is taken from the set of eigenvalues of the operator $A_{0}$ that are contained in the set $(-\infty,\theta_{k})$. We intend to prove the following theorems that provide sufficient conditions for the existence of compact full solutions connecting stationary points for the system \eqref{sys-res} in the case of the resonance at infinity.
\begin{theorem}\label{th-crit-ogr}
Let us assume that $\{h_{k}\}_{k=1}^{m}$ is a family of $L^{2}(\Omega)$ functions such that the inequalities $(C1)_{\pm}$ and $(C2)_{\pm}$ are fulfilled and suppose that the following non-resonance condition at the origin holds
\begin{equation*}
\sigma(A_{0})\cap\sigma(G+\Lambda)=\emptyset.
\end{equation*}
If the resonance conditions $(LL1)_{\pm}$ and $(LL2)_{\pm}$ are satisfied, then there is a non-trivial bounded full solution $u$ of the system \eqref{sys-res} such that either $u(t)\to 0$ as $t\to+\infty$ or $u(t)\to 0$ as $t\to-\infty$, provided 
\begin{equation}\label{cp-1}
d_{0}(\lambda)\neq d_{\infty}(\lambda)+(n_{1}(\lambda)+n_{2}(\lambda))/2\pm (n_{1}(\lambda)+ n_{2}(\lambda))/2.
\end{equation}
\end{theorem}
\begin{theorem}\label{th-crit-ogr-2}
Let us assume that $\{h_{k}\}_{k=1}^{m}$ is a family of $L^{2}(\Omega)$ functions such that the inequalities $(C1)_{\pm}$ and $(C2)_{\mp}$ are fulfilled and suppose that the following non-resonance condition at the origin holds
\begin{equation}\label{non-res-origin}
\sigma(A_{0})\cap\sigma(G+\Lambda)=\emptyset.
\end{equation}
If the resonance conditions $(LL1)_{\pm}$ and $(LL2)_{\mp}$ are satisfied, then there is a non-trivial bounded full solution $u$ of the system \eqref{sys-res} such that either $u(t)\to 0$ as $t\to+\infty$ or $u(t)\to 0$ as $t\to-\infty$, provided 
\begin{equation}\label{cp-2}
d_{0}(\lambda)\neq d_{\infty}(\lambda)+(n_{1}(\lambda)+n_{2}(\lambda))/2\pm (n_{1}(\lambda)- n_{2}(\lambda))/2.
\end{equation}
\end{theorem}
\noindent{\em Proof of Theorems \ref{th-crit-ogr} and \ref{th-crit-ogr-2}.} Observe that the operator $L$, given by the formula
$$L := (A_{0} - \lambda_{1}I)\times\ldots\times (A_{0} - \lambda_{m}I)- DF'(0),$$
is sectorial and has compact resolvents as \cite[Proposition 3.1.4]{Pazy} and \cite[Theorem 3.2.1]{Pazy} say. Then the spectrum of $L$ consists of a sequence of real eigenvalues, which is either finite or diverges to infinity. Furthermore the respective eigenspaces are finite dimensional. We claim that 
\begin{align}\label{numb-d-0}
d_{0}(\lambda)= \sum_{\mu<0}\dim\ker(\mu I - L),
\end{align}
where, in the above summation, the parameter $\mu$ is taken from the set of all negative eigenvalues of the operator $L$. Indeed, in view of \eqref{der-11}, we have
\begin{align*}
L u = (A_{0} u_{1},\ldots,A_{0} u_{m}) -(G+\Lambda) u,\quad  u\in D(L).
\end{align*}
If we take $O$ to be the $m\times m$ orthogonal matrix such that $$O^{t}(G+\Lambda)O = \mathrm{diag}\,(\theta_{1},\theta_{2},\ldots,\theta_{m}),$$ then, for any $ u\in D(L)$, the following holds
\begin{equation}
\begin{aligned}\label{ker-11}
L(O u) &= (A_{0}(O u)_{1},\ldots, A_{0} (O u)_{m}) -(G+\Lambda)(O u) \\
&= O(A_{0}  u_{1},\ldots,A_{0}  u_{m}) - O(\theta_{1}  u_{1},\ldots, \theta_{m}  u_{m}),\\
&= O[(A_{0}-\theta_{1}I) u_{1},\ldots,(A_{0}-\theta_{m}I) u_{m}],
\end{aligned}
\end{equation}
which implies that $\sigma(L) = \sigma((A_{0} - \theta_{1} I)\times\ldots\times(A_{0} - \theta_{m} I))$. Consequently 
\begin{align}\label{eq-kk-pp}
\sum_{\mu<0}\dim\ker(\mu I - L) = \sum_{\mu<0} \sum_{k=1}^{m}\dim\ker((\mu+\theta_{k})I-A_{0}),
\end{align}
where in the above summations the parameters $\mu$ ranges over the negative elements of the set $\sigma(L)$. Let us observe that, for any $1\le k\le m$, we have
\begin{gather*}
\{\mu + \theta_{k} \ | \ \mu\in\sigma(L), \ \mu<0 \text{ and }\mathrm{Ker}\,((\mu + \theta_{k})I - A_{0}) \neq \{0\}\} \\
= \{\nu<\theta_{k} \ | \ \mathrm{Ker}\,(\nu I - A_{0})\neq\{0\}\}.
\end{gather*}
Combining this with \eqref{eq-kk-pp}, we obtain
\begin{align*}
&\sum_{\mu<0}\dim\ker(\mu I - L)= \sum_{k=1}^{m}\sum_{\mu<0}\dim\ker((\mu+\theta_{k})I-A_{0}) \\
&\qquad\quad= \sum_{k=1}^{m} \sum_{\nu<\theta_{k}}\dim\ker(\nu I - A_{0})=d_{0}(\lambda),
\end{align*}
where in the above summations the parameters $\mu$ and $\nu$ ranges over the sets $\sigma(L)$ and $\sigma(A_{0})$, respectively. This gives \eqref{numb-d-0} as desired. Let us observe that the condition \eqref{non-res-origin} implies that $\ker[(A_{0}-\theta_{1}I)\times\ldots\times(A_{0}-\theta_{m}I)] = \{0\}$, which together with \eqref{ker-11} yield $\ker L = \{0\}$. Therefore, from \cite[Theorem II.3.5]{MR910097} it follows that the invariant set $K_{0}:=\{0\}$ admits an admissible isolating neighborhood and $h(\Phi,K_0)= \Sigma^{d_{0}(\lambda)}$, where $\Phi$ is the semiflow associated with the system \eqref{sys-res}. If the inequalities $(C1)_{\pm}$, $(C2)_{\pm}$ and conditions $(LL1)_{\pm}$, $(LL2)_{\pm}$ are satisfied, then, by Theorem \ref{th-1}, we infer that the set $K_{\infty}$ consisting of all bounded full solutions of the semiflow 
$\Phi$ also has an admissible isolating neighborhood and its homotopy index is given by 
\begin{equation*}
h(\Phi, K_{\infty}) = \left\{\begin{aligned}
&\Sigma^{d_{0}(\lambda)+n_{1}(\lambda)+n_{2}(\lambda)}&& \!\text{if $(C1)_{+}$, $(C2)_{+}$, $(LL1)_{+}$, $(LL2)_{+}$ hold},\\
&\Sigma^{d_{0}(\lambda)}&& \!\text{if $(C1)_{-}$, $(C2)_{-}$, $(LL1)_{-}$, $(LL2)_{-}$ hold}.
\end{aligned}\right.
\end{equation*}
Clearly $K_0\subset K_{\infty}$ and $h(\Phi, K_{0})\neq \o 0$. By the condition \eqref{cp-1}, we have also that $h(\Phi, K_{0})\neq h(\Phi, K_{\infty})$. Hence Proposition \ref{prop-irred} gives the existence of a non-trivial full solution $u$ of the semiflow $\Phi$ such that either $u(t)\to 0$ as $t\to+\infty$ or $u(t)\to 0$ as $t\to-\infty$. Thus the proof of Theorems \ref{th-crit-ogr} is completed. On the other hand, if the inequalities $(C1)_{\pm}$, $(C2)_{\mp}$ together with the conditions $(LL1)_{\pm}$, $(LL2)_{\mp}$ are satisfied, then Theorem \ref{th-2} says that the set $K_{\infty}$ admits an admissible isolating neighborhood and its homotopy index is such that
\begin{equation*}
h(\Phi, K_{\infty}) = \left\{\begin{aligned}
&\Sigma^{d_{0}(\lambda)+n_{1}(\lambda)}&& \!\text{if $(C1)_{+}$, $(C2)_{-}$, $(LL1)_{+}$, $(LL2)_{-}$ hold},\\
&\Sigma^{d_{0}(\lambda)+n_{2}(\lambda)}&& \!\text{if $(C1)_{-}$, $(C2)_{+}$, $(LL1)_{-}$, $(LL2)_{+}$ hold}.
\end{aligned}\right.
\end{equation*}
Consequently the condition \eqref{cp-2}, implies that again $h(\Phi, K_{0})\neq h(\Phi, K_{\infty})$ and therefore Proposition \ref{prop-irred} provides the existence of a non-trivial full solution $u$ of the semiflow $\Phi$ such that either $u(t)\to 0$ as $t\to+\infty$ or $u(t)\to 0$ as $t\to-\infty$. This in turn completes the proof of Theorem \ref{th-crit-ogr-2}. \hfill $\square$ \\

In Theorems \ref{th-crit-ogr} and \ref{th-crit-ogr-2} we employed the resonance conditions $(LL1)_{\pm}$ and $(LL2)_{\pm}$ to derive the existence of a full solution $u$ of the semiflow $\Phi$ with the relatively compact image $u(\R)\subset X^{\alpha}$ such that $0\in \alpha(u)\cup\omega(u)$. To deduce that $u$ connects the origin with a nontrivial stationary point we make an additional assumption on the semiflow $\Phi$. 
\begin{definition}
We say that the semiflow $\Phi$ is {\em gradient-like} with respect to the functional $V:X^{\alpha}\to\R$ provided $$V(\Phi(u_{0},t_{1}))\ge V(\Phi(u_{0},t_{2})),\quad u_{0}\in X^{\alpha}, \ t_{1}>t_{2}\ge 0$$ and, for any non-constant full solution $u$ of $\Phi$, the value $V(u(t))$ is not constant for $t\in\R$. Then $V$ is called the {\em Liapunov function} for the semiflow $\Phi$. \hfill $\square$
\end{definition}
\begin{remark}
The usual assumption which makes $\Phi$ a gradient-like semiflow is the existence of a smooth potential function $\tilde f:\Omega\times\R^{m}\to\R$ such that
\begin{align*}
f_{k}(x,s,y) = \partial_{s_{k}}\tilde f(x,s), \quad x\in\Omega, \ s\in\R^{m}, \ y\in\R^{nm}, \ 1\le k\le m.
\end{align*}
Then the energy functional $E:X^{\alpha}\to\R$ given, for any $u=(u_{1},\ldots,u_{m})\in X^{\alpha}$, by 
\begin{align*}
E(u):=\sum_{k=1}^{m} \frac{1}{2}\int_{\Omega}|\nabla u_{k}(x)|^{2}\,dx + \int_{\Omega} \tilde f(u(x))\,dx
\end{align*}
is the Liapunov function for the semiflow $\Phi$ (see \cite{MR1778284} for more details). \hfill $\square$
\end{remark}
The following corollaries are simple consequences of Theorems \ref{th-crit-ogr} and \ref{th-crit-ogr-2}.
\begin{corollary}\label{th-crit-ogr3}
Let us assume that the semiflow $\Phi$ is gradient-like and 
\begin{equation*}
\sigma(A_{0})\cap\sigma(G+\Lambda)=\emptyset.
\end{equation*}
Suppose that $\{h_{k}\}_{k=1}^{m}$ is a family of $L^{2}(\Omega)$ functions satisfying the inequalities $(C1)_{\pm}$ and $(C2)_{\pm}$. If the resonance conditions $(LL1)_{\pm}$ and $(LL2)_{\pm}$ hold and 
\begin{equation*}
d_{0}(\lambda)\neq d_{\infty}(\lambda)+(n_{1}(\lambda)+n_{2}(\lambda))/2\pm (n_{1}(\lambda)+ n_{2}(\lambda))/2,
\end{equation*}
then there are non-zero stationary point $u_{0}\in X^{\alpha}$ and full solution $u$ of the system \eqref{sys-res} such that either $u(t_{n}) \to  u_{0}$ for some $t_{n}\to+\infty$ and $u(t)\to 0$ as $t\to-\infty$ or $u(t_{n}) \to  u_{0}$ for some $t_{n}\to-\infty$ and $u(t)\to 0$ as $t\to+\infty$.
\end{corollary}
\begin{corollary}\label{th-crit-ogr4}
Let us assume that the semiflow $\Phi$ is gradient-like and 
\begin{equation*}
\sigma(A_{0})\cap\sigma(G+\Lambda)=\emptyset.
\end{equation*}
Suppose that $\{h_{k}\}_{k=1}^{m}$ is a family of $L^{2}(\Omega)$ functions satisfying the inequalities $(C1)_{\pm}$ and $(C2)_{\mp}$. If the resonance conditions $(LL1)_{\pm}$ and $(LL2)_{\mp}$ hold and 
\begin{equation*}
d_{0}(\lambda)\neq d_{\infty}(\lambda)+(n_{1}(\lambda)+n_{2}(\lambda))/2\pm (n_{1}(\lambda)- n_{2}(\lambda))/2,
\end{equation*}
then there are non-zero stationary point $u_{0}\in X^{\alpha}$ and full solution $u$ of the system \eqref{sys-res} such that either $u(t_{n}) \to u_{0}$ for some $t_{n}\to+\infty$ and $u(t)\to 0$ as $t\to-\infty$ or $u(t_{n}) \to u_{0}$ for some $t_{n}\to-\infty$ and $u(t)\to 0$ as $t\to+\infty$.
\end{corollary}
\noindent{\em Proof of Corollaries \ref{th-crit-ogr3} and \ref{th-crit-ogr4}.} In view of Theorems \ref{th-crit-ogr} and \ref{th-crit-ogr-2}, we infer that there is a non-trivial full solution $u$ of the system \eqref{sys-res} such that either $u(t)\to 0$ as $t\to+\infty$ or $u(t)\to 0$ as $t\to-\infty$. Since the semiflow $\Phi$ is gradient-like, from \cite[Theorem II.5.4]{MR910097} we know that the the limit sets $\alpha(u)$ and $\omega(u)$ are non-empty, disjoint and consist of the stationary points of the semiflow $\Phi$. Consequently, if $u(t)\to 0$ as $t\to+\infty$, then we can choose a non-zero element $u_{0}\in\alpha(u)$ such that $u(t_{k})\to  u_{0}$ for some $t_{k}\to -\infty$. On the other hand, if $u(t)\to 0$ as $t\to-\infty$, then we can take $u_{0}\in\omega(u)$ such that $u_{0}\neq 0$ and $u(t_{k})\to u_{0}$ for some $t_{k}\to+\infty$. Thus the proof of the corollaries is complete. \hfill $\square$

\def\cprime{$'$} \def\polhk#1{\setbox0=\hbox{#1}{\ooalign{\hidewidth
  \lower1.5ex\hbox{`}\hidewidth\crcr\unhbox0}}} \def\cprime{$'$}
  \def\cprime{$'$} \def\cprime{$'$}
\providecommand{\bysame}{\leavevmode\hbox to3em{\hrulefill}\thinspace}
\providecommand{\MR}{\relax\ifhmode\unskip\space\fi MR }
\providecommand{\MRhref}[2]{%
  \href{http://www.ams.org/mathscinet-getitem?mr=#1}{#2}
}
\providecommand{\href}[2]{#2}

\parindent = 0 pt


\begin{thebibliography}{10}

\bibitem{MR1399539}  A. Ambrosetti, D. Arcoya, {\em On a quasilinear problem at strong resonance}, Topol. Methods Nonlinear Anal.  {\bf 6}  (1995),  no. 2, 255--264.

\bibitem{MR0492839} A. Ambrosetti, G. Mancini, {\em Existence and multiplicity results for nonlinear elliptic problems with linear part at resonance. The case of the simple eigenvalue},  J. Differential Equations  {\bf 28}  (1978),  no. 2, 220--245.

\bibitem{MR0487001} A. Ambrosetti, G. Mancini, \emph{Theorems of existence and multiplicity for nonlinear elliptic problems with noninvertible linear part}, Ann. Scuola Norm. Sup. Pisa Cl. Sci. (4) \textbf{5} (1978), no.~1, 15--28.

\bibitem{MR2280977} P. Amster, P. De N\'apoli, Pablo . {\em Landesman-Lazer type conditions for a system of $p$-Laplacian like operators}, J. Math. Anal. Appl. 326 (2007), no. 2, 1236--1243.

\bibitem{MR1041504} D. Arcoya, A. Ca\~{n}ada, {\em Critical point theorems and applications to nonlinear boundary value problems}, Nonlinear Anal. 14  (1990),  no. 5, 393--411.

\bibitem{MR1296115} D. Arcoya, D.G. Costa, {\em Nontrivial solutions for a strongly resonant problem}, Differential Integral Equations  8  (1995),  no. 1, 151--159.

\bibitem{MR1430505} D. Arcoya, L. Orsina, {\em Landesman-Lazer conditions and quasilinear elliptic equations},
 Nonlinear Anal. {\bf 28}  (1997),  no. 10, 1623--1632.

\bibitem{bredon} G.E. Bredon, {\em Topology and Geometry}, Springer-Verlag, 1993.

\bibitem{MR1219185} D.G. Costa,  E.A. Silva, {\em Existence of solutions for a class of resonant elliptic problems}, J. Math. Anal. Appl.  175  (1993),  no. 2, 411--424.

\bibitem{arieta} J. Arrieta, R. Pardo, A. Rodriguez-Bernal, \emph{Equilibria and global dynamics of a problem with bifurcation from infinity}, J. Differential Equations 246 (2009), 2055-2080.

\bibitem{MR1756573}  J. Bouchala, P. Dr\'abek, {\em Strong resonance for some quasilinear elliptic equations},
 J. Math. Anal. Appl.  {\bf 245}  (2000),  no. 1, 7--19.

\bibitem{MR713209} P. Bartolo, V. Benci, D. Fortunato, \emph{Abstract critical point theorems and applications to some nonlinear problems with ``strong'' resonance at infinity}, Nonlinear Anal. \textbf{7} (1983), no.~9, 981--1012.

\bibitem{MR0513090} H. Br{\'e}zis, L. Nirenberg, \emph{Characterizations of the ranges of some nonlinear operators and applications to boundary value problems}, Ann. Scuola Norm. Sup. Pisa Cl. Sci. (4) \textbf{5} (1978), no.~2, 225--326.

\bibitem{MR0448180} L. Cesari, R. Kannan, {\em An abstract existence theorem at resonance}, Proc. Amer. Math. Soc.  63  (1977),  no. 2, 221--225.

\bibitem{MR1778284}
J.~W. Cholewa, T. D{\l}otko, \emph{Global attractors in abstract parabolic problems}, London Mathematical Society Lecture Note Series, vol. 278,
  Cambridge University Press, Cambridge, 2000.

\bibitem{MR511133}
C. Conley, \emph{Isolated invariant sets and the {M}orse index}, CBMS
  Regional Conference Series in Mathematics, vol.~38, American Mathematical
  Society, Providence, R.I., 1978.

\bibitem{MR3883315} A. \'Cwiszewski, W. Kryszewski, {\em Bifurcation from infinity for elliptic problems on $\R^{N}$}, Calc. Var. Partial Differential Equations  {\bf 58}  (2019),  no. 1, Art. 13,

\bibitem{MR1151266} D.G. de~Figueiredo, J.P. Gossez, \emph{Strict monotonicity of eigenvalues and unique continuation}, Comm. Partial Differential Equations \textbf{17} (1992), no.~1-2, 339--346.

\bibitem{MR2085538} P. Dr\'abek, P. Girg, P. Tak\'a\v{c}, {\em Bounded perturbations of homogeneous quasilinear operators using bifurcations from infinity}, J. Differential Equations  204  (2004),  no. 2, 265--291.

\bibitem{MR1726752} P. Dr\'abek, S.B. Robinson, {\em Resonance problems for the $p$-Laplacian}, J. Funct. Anal. {\bf 169}  (1999),  no. 1, 189--200.

\bibitem{MR1689320} P. Dr\'abek, S.B Robinson, {\em Resonance problems for the one-dimensional $p$-Laplacian}, Proc. Amer. Math. Soc.  128  (2000),  no. 3, 755--765.

\bibitem{MR0541873} D.G. de Figueiredo, W.M. Ni, {\em Perturbations of second order linear elliptic problems by nonlinearities without Landesman-Lazer condition}, Nonlinear Anal.  3  (1979),  no. 5, 629--634.

\bibitem{MR2609541}  P. Felmer, A. Quaas, B. Sirakov, {\em Landesman-Lazer type results for second order Hamilton-Jacobi-Bellman equations}, J. Funct. Anal.  {\bf 258}  (2010),  no. 12, 4154--4182.

\bibitem{MR2810138} A. Fonda, M. Garrione, {\em Nonlinear resonance: a comparison between Landesman-Lazer and Ahmad-Lazer-Paul conditions},
 Adv. Nonlinear Stud.  {\bf 11}  (2011),  no. 2, 391--404.

\bibitem{MR882069} N. Garofalo, F.H. Lin, \emph{Unique continuation for elliptic operators: a geometric-variational approach}, Comm. Pure Appl. Math. \textbf{40} (1987), no.~3, 347--366.

\bibitem{hat} A. Hatcher, {\em Algebraic Topology}, Cambridge Univ. Press, 2002.

\bibitem{MR610244} D. Henry, \emph{Geometric theory of semilinear parabolic equations}, Lecture Notes in Mathematics, vol. 840, Springer-Verlag, Berlin, 1981.

\bibitem{H-F} E. Hille, R. Phillips, {\em Functional Analysis and Semi-Groups}, Colloquium Publications 31, American Mathematical Society, Providence, RI, 1957.

\bibitem{MR1466317} J. Johnsen, T. Runst, {\em Semi-linear boundary problems of composition type in $L^{p}$-related spaces}, Comm. Partial Differential Equations  22  (1997),  no. 7-8, 1283--1324.

\bibitem{Kok5} P. Kokocki, \emph{Connecting orbits for nonlinear differential equations at resonance.} J. Differential Equations 255 (2013), no. 7, 1554-1575

\bibitem{MR0267269}
E.M. Landesman, A.C. Lazer, \emph{Nonlinear perturbations of linear elliptic boundary value problems at resonance}, J. Math. Mech. \textbf{19}   (1969/1970), 609--623.

\bibitem{MR1621026}  S. Ma, Z. Wang, J. Yu, {\em An abstract existence theorem at resonance and its applications}, J. Differential Equations  {\em 145}(1998),  no. 2, 274--294.

\bibitem{MR0956010} J. Mawhin, K. Schmitt, {\em Landesman-Lazer type problems at an eigenvalue of odd multiplicity},
 Results Math.  {\bf 14}  (1988),  no. 1-2, 138--146.

\bibitem{MR1306589} R. Ortega, {\em A boundedness result of Landesman-Lazer type}, Differential Integral Equations  8  (1995),  no. 4, 729--734.

\bibitem{MR1342038} R. Ortega, A. Tineo, \emph{Resonance and non-resonance in a problem of boundedness}, Proc. Amer. Math. Soc. \textbf{124} (1996), no.~7,  2089--2096.

\bibitem{Pazy} A. Pazy, {\em Semigroups of linear operators and applications to partial differential equations}, Springer Verlag 1983.

\bibitem{MR1992823} M. Prizzi, \emph{On admissibility for parabolic equations in {$\Bbb R\sp n$}}, Fund. Math. \textbf{176} (2003), no.~3, 261--275.

\bibitem{MR637695}
K.P. Rybakowski, \emph{On the homotopy index for infinite-dimensional semiflows}, Trans. Amer. Math. Soc. \textbf{269} (1982), no.~2, 351--382.

\bibitem{MR798176}
K.P. Rybakowski, \emph{Nontrivial solutions of elliptic boundary value problems with resonance at zero}, Ann. Mat. Pura Appl. (4) \textbf{139} (1985), 237--277.

\bibitem{MR910097}
K.P. Rybakowski, \emph{The homotopy index and partial differential equations}, Universitext, Springer-Verlag, Berlin, 1987.

\bibitem{MR731150}
K.P. Rybakowski, \emph{Trajectories joining critical points of nonlinear parabolic and hyperbolic partial differential equations}, J. Differential Equations \textbf{51} (1984), no.~2, 182--212.

\bibitem{MR797044}
D. Salamon, \emph{Connected simple systems and the {C}onley index of isolated invariant sets}, Trans. Amer. Math. Soc. \textbf{291} (1985), no.~1,
  1--41.

\bibitem{MR1055653} M. Schechter, {\em Nonlinear elliptic boundary value problems at strong resonance}, Amer. J. Math. {\bf 112} (1990), no. 3, 439--460.

\bibitem{MR688146} J. Smoller, \emph{Shock waves and reaction-diffusion equations}, Grundlehren der Mathematischen Wissenschaften, vol. 258, Springer-Verlag, New York, 1983.

\bibitem{spanier} E.H. Spanier, {\em Algebraic topology}, Mc-Graw-Hill, New York, 1996.

\bibitem{MR0566069} K. Thews, {\em Nontrivial solutions of elliptic equations at resonance}, Proc. Roy. Soc. Edinburgh Sect. A  85  (1980),  no. 1-2, 119--129.

\bibitem{MR0612616} K. Thews, {\em $T$-periodic solutions of time dependent Hamiltonian systems with a potential vanishing at infinity}, Manuscripta Math. 33 (1980/81), no. 3-4, 327--338.

\bibitem{MR500580}
H. Triebel, \emph{Interpolation theory, function spaces, differential operators}, VEB Deutscher Verlag der Wissenschaften, Berlin, 1978.

\bibitem{MR597281}
J. Valdo, A. Gon{\c{c}}alves, \emph{On bounded nonlinear perturbations of an elliptic equation at resonance}, Nonlinear Anal. \textbf{5} (1981), no.~1,
  57--60.

\end{thebibliography}
\end{document}